\documentclass[12pt]{mcom-l}
\newtheorem{thm}{Theorem}[section]
\newtheorem{cor}{Corollary}[section]
\newtheorem{lem}{Lemma}[section]

\newtheorem{conj}{Conjecture}[section]
\newtheorem{prop}{Proposition}[section]
\newtheorem{cexample}{Counterexample}[section]

\newcommand{\BbN}{\mathbb{N}}
\newcommand{\BbZ}{\mathbb{Z}}
\newcommand{\BbQ}{\mathbb{Q}}
\newcommand{\Res}{\mathrm{Res}}
\newcommand{\F}{{\mathcal F}}
\newcommand{\A}{{\mathcal A}}
\newcommand{\tm}{t_{\mathrm{M}}}
\newcommand{\tz}{t_{\BbZ}}
\newcommand{\tr}{\mathrm{cap}}
\newcommand{\Lm}{L_{-}}
\newcommand{\Lp}{L_{+}}
\newcommand{\al}{\alpha}
\newcommand{\lh }{\|}
\newcommand{\rh}{\|^{{}^{\scriptstyle{*}}}}
\newcommand{\Farey}{\left[\tfrac{b_1}{c_1},\tfrac{b_2}{c_2}\right]}
\newcommand{\fracn}{\frac{2}{n+\sqrt{n^2-4}}}
\usepackage{psfig}
\usepackage{amsfonts}
\usepackage{amssymb}
\usepackage{enumerate}
\usepackage{comment}

\title[Monic integer transfinite diameter]{The monic integer transfinite diameter}
\author{K. G. Hare}
\address{Department of Pure Mathematics, University of Waterloo, Waterloo,
    Ontario, Canada,  N2L 3G1}
\email{kghare@math.uwaterloo.ca}
\thanks{Research of K.G. Hare was supported in part by NSERC of Canada
    and a Seggie Brown Fellowship, University of Edinburgh.}

\author{C. J. Smyth}
\address{School of Mathematics, University of Edinburgh,
    James Clerk Maxwell Building, King's Buildings,
    Mayfield Road, Edinburgh EH9 3JZ, UK.}
\email{c.smyth@ed.ac.uk}

\begin{document}

\begin{abstract}

We study the problem of finding nonconstant monic integer
polynomials, normalized by their degree, with
    small  supremum   on an interval $I$. The monic integer transfinite diameter $\tm(I)$ is defined as
     the infimum of all such supremums.
We show that if $I$ has length $1$ then $\tm(I) = \tfrac{1}{2}$.

We make three general conjectures relating to the value of $\tm(I)$
for
    intervals $I$ of length less that $4$. We also conjecture a value for $\tm([0, b])$ where $0<b\le 1$.
We  give some partial
    results, as well as computational evidence, to support these conjectures.

We define functions $\Lm(t)$ and $\Lp(t)$, which measure properties
of the lengths
    of intervals $I$ with $\tm(I)$ on either side of $t$.
Upper and lower bounds are given for these functions.

We also consider the problem of determining $\tm(I)$ when
    $I$ is a Farey interval.
We prove that a conjecture of Borwein,
    Pinner and Pritsker concerning this value is true for an infinite family of Farey intervals.

\end{abstract}



\maketitle
\section{Introduction and Results}

In this paper we continue a study, recently initiated  by Borwein, Pinner and
    Pritsker \cite{BorweinPinnerPritsker03}, of the
    {\em monic integer transfinite diameter} of a real interval.
We write the normalized supremum  on an interval $I$ as
    \[ \lh P\rh_I := \sup_{x \in I} |P(x)|^{1/\deg P}. \]
     Note that this is not a norm. Then the monic integer transfinite diameter $\tm(I)$ is defined as
    \[ \tm(I) := \inf_P \lh P\rh_I,\]
    where the infimum is taken over all
    non-constant monic polynomials with integer coefficients. We call $\tm(I)$ the {\it monic integer transfinite diameter} of $I$ (also called the {\it monic integer Chebyshev constant} \cite{Borwein02,BorweinPinnerPritsker03}).
Clearly $\tm(I)\ge \tz(I)$, where $\tz(I)$ denotes the
    {\em  integer transfinite diameter}, defined using the same infimum,
    but taken over the larger set of all non-constant polynomials with
    integer coefficients \cite{BorweinErdelyi96,Chudnovsky83,FlammangRhinSmyth97}.
Further $\tz(I) \ge \tr(I)$, the {\em capacity} or
    {\em transfinite diameter} of $I$ \cite{Goluzin69, Ransford95},
    which can be defined again using the same
    infimum, but this time taken over all non-constant monic polynomials
    with real coefficients.
It is well known that $\tr(I)=|I|/4$ for an interval $I$ of length $|I|$.
Further, if $|I|\ge 4$ then $\tz(I)=\tm(I)=\tr(I)$ by
    \cite{BorweinPinnerPritsker03}
    so that the challenge for evaluating $\tm(I)$, as for $\tz(I)$, lies in
     intervals with $|I|<4$. For these intervals we know from
     \cite[Prop. 1.2]{BorweinPinnerPritsker03} that $\tm(I)<1$.
However, in contrast to the study of $\tz(I)$,
    in the monic case it is possible to evaluate $\tm(I)$ exactly over
    some such intervals.

Our first result is the following.
\begin{thm}
\label{thm:interval 1} All intervals $I$ of length $1$ have $\tm(I)
= \tfrac{1}{2}$. In fact, slightly more is true:
    if $1 \leq |I| \leq 1.008848$ then $\tm(I) = \tfrac{1}{2}$.

    Furthermore
    for any $b<1$ there is an interval $I$ with $|I|=b$ and
    $\tm(I) < \tfrac{1}{2}$, while for $b>1.064961507$
    there is an interval $I$ with $|I|=b$
    and $\tm(I) > \tfrac{1}{2}$.
\end{thm}

The proof, which is essentially a corollary of Theorem
 \ref{thm:\Lm(1/2)} (a) below, is discussed in Section \ref{sec:interval 1}.

The numbers, 1.008848 and 1.064961507 in Theorem \ref{thm:interval 1},
    like most numerical values given
    in this paper, are approximations to some exact algebraic number.
These numbers are rounded in the correct direction, if necessary, to
    ensure an inequality still holds.
 The  polynomial equations that they  satisfy is given within the text.
We  have tried to do this for all numerical values.

    To measure the range of lengths of intervals having a particular monic integer transfinite diameter $t$,
    we introduce the following two functions:
\begin{eqnarray*}
\Lm(t)&:=&\inf_I\{|I|: \tm(I)>t\};\\
\Lp(t)&:=&\sup_I\{|I|: \tm(I)\leq t\}.
\end{eqnarray*}

It follows from \cite[Prop. 1.3]{BorweinPinnerPritsker03} that both
$\Lm(t)$ and $\Lp(t)$ are nondecreasing functions of $t$. Also
$\Lm(t)\le\Lp(t)$ -- see Lemma \ref{lem:LL}(a) below.
 We give (Proposition \ref{prop:L-bounds})
general method for finding upper and lower bounds for $\Lm(t)$ and
$\Lp(t)$, and apply these methods to get such bounds for $\tfrac{1}{2} \leq t
\leq 1$. They are constructive, using both the LLL basis-reduction
algorithm
 and the Simplex method. These techniques were first applied in this area by Borwein and
Erd\'elyi \cite{BorweinErdelyi96}, and then
    by Habsieger and Salvy \cite{HabsiegerSalvy97}.
   These bounds are given in Theorem \ref{thm:\Lm(t)} and Proposition \ref{P-lowerB} -- see also
   Figures \ref{fig:\Lm(t)} and \ref{fig:\Lp(t)}.

At $t=\tfrac{1}{2}$, we pushed this method further, and were able to say
more.

\begin{thm}
\label{thm:\Lm(1/2)} We have
\begin{itemize}
\item[(a)]
 $1.008848 \leq \Lm\left(\tfrac{1}{2}\right) \leq 1.064961507$

\noindent and
\item[(b)] $\sqrt{2} \approx 1.41421 \leq
\Lp\left(\tfrac{1}{2}\right) \leq 1.4715$.
\end{itemize}
\end{thm}

Further properties of $\Lp$ and $\Lm$ are given in Lemma
\ref{lem:LL}.


\section{Definitions, Conjectures and Further Results}

In this section, we state some old and some more new results, and
    (perhaps a little recklessly) make four conjectures.

The following result is simple but fundamental. It is
      useful for determining lower bounds for $\tm(I)$.

\vspace{0.1 in}

 \noindent {\bf Lemma BPP} (Borwein, Pinner and Pritsker
   \cite[p.1905]{BorweinPinnerPritsker03}).
\begin{em}
Let $Q(x) = a_d x^d + \cdots + a_0$ be a nonmonic irreducible
polynomial with
    integer coefficients, all of whose roots lie in
    the interval $I$.
Then $\lh P\rh_I \geq a_d^{-1/d}$ for every monic integer polynomial $P$, so that $\tm(I)\geq a_d^{-1/d}$.
Furthermore, if $\lh P\rh_I = a_d^{-1/d}$ then $\tm(I) = a_d^{-1/d}$
and
    $|P(\beta)|^{1/\deg P}=a_d^{-1/d}$ for every root $\beta$ of
    $Q$, and $\Res(P,Q)=\pm 1$.
\end{em}

\vspace{0.1 in}

The proof follows straight from the classical fact that, for the
conjugates $\beta_i$ of $\beta$
\begin{equation}\label{eqn:resultant}
\Res(P,Q)=a_d^{\deg P}\prod_{i=1}^d P(\beta_i)
\end{equation}
is a nonzero integer, giving
\begin{equation}\label{eqn:resultant2}
\lh P\rh_I\ge \left(\prod_i|P(\beta_i)|^{1/\deg
P}\right)^{\frac{1}{d}}\ge a_d^{-1/d}|\Res(P,Q)|^{\frac{1}{d\deg
P}}\ge a_d^{-1/d}.
\end{equation}
This result is a variant of a similar one in the theory of
$\tz(I)$---see Lemma \ref{lem:was BE}.

We call such a value $a_d^{-1/d}$ in Lemma BPP an {\em obstruction}
for $I$, with
 {\em obstruction polynomial} $Q(x)$.  From
Lemma BPP we see that $\tm(I)$ is bounded below by the supremum
    of all such obstructions.
If this supremum is attained by some value $a_d^{-1/d}$ coming from
$Q(x)=a_dx^d+\cdots+a_0$, then we say $a_d^{-1/d}$ is a {\em maximal
obstruction}, and $Q(x)$ is a {\em maximal obstruction
    polynomial}.
It is not known whether such a polynomial exists for all intervals
$I$ of length less than $4$ (see
    Conjecture \ref{conj:maximal}).

    We say that  the monic integer polynomial $P(x)$ is an {\em optimal monic
    integer Chebyshev polynomial for $I$}  if  $\lh P\rh_I=\tm(I)$.
If $I$ has a maximal obstruction  $a_d^{-1/d}$  with
$\tm(I)=a_d^{-1/d}$
     and an optimal monic integer Chebyshev polynomial $P$ then we say that
    {\em $P$ attains the maximal obstruction $a_d^{-1/d}$}.

Throughout this paper, $P(x)$ will denote a monic integer
    polynomial,  $Q(x)$ a nonmonic integer polynomial and $R(x)$ any integer
    polynomial.

One very nice property of the monic integer transfinite diameter problem,
    not shared by its nonmonic cousin, is that often
    exact values can be computed for $\tm(I)$.
In all cases where this has been done, including Theorem \ref{thm:interval 1}, it was achieved
    by finding a maximal obstruction, and a corresponding optimal monic integer Chebyshev polynomial.
    Simple examples of this are given (\cite[Theorem 1.5]{BorweinPinnerPritsker03}) by the intervals $I=[0,1/n]$  for
$n\ge 2$, where $Q(x)=nx-1$ is a maximal obstruction polynomial, and
$P(x)=x$ is an optimal monic integer Chebyshev polynomial. For
$n=1$, $\tm([0,1])=\tfrac{1}{2}$, with $Q(x)=2x-1$ and
$P(x)=x(x-1)$. This was  the case too in \cite[Section
5]{BorweinPinnerPritsker03} in the proof
    of the Farey Interval conjecture for small-denominator intervals.

A much less obvious example is the interval $I = [-0.3319, 0.7412]$,
    of length $1.0731$.
Here, we have  $\tm(I)=\lh P\rh_I= 7^{-1/3} \approx 0.522$,
     with maximal obstruction polynomial
    $7 x^3 - 7 x^2 + 1$
and where $P$ is the optimal monic integer Chebyshev polynomial
\[
\begin{array}{rl}
P(x)=&x^{276507}(x-1)^{29858} (x^2+x-1)^{14929}\\
 &(x^5-17 x^4+24 x^3-8 x^2-2 x+1)^{28848} \\
&(x^7-117 x^6+194 x^5-70 x^4-31 x^3+18 x^2+x-1)^{7935} \\
&(x^8-4 x^7+97 x^6-172 x^5+78 x^4+20 x^3-18 x^2+1)^{9795} \\
&(x^8-34 x^7+164 x^6-208 x^5+65 x^4+33 x^3-18 x^2-x+1)^{5846} \\
&(x^8-7 x^7+2 x^6-x^5-10 x^4+28 x^3-15 x^2-2 x+2)^{1148}
\end{array}
\]
    of degree $670320$.
(Tighter endpoints for this interval, and its length, can be computed
by
    solving the equation $P(x) = \pm \left(7^{-1/3}\right)^{\deg P}$.)
The discovery of this polynomial required the use of
    Lemma \ref{lem:resultant-value} below.

For the nonmonic transfinite diameter $\tz$, Pritsker \cite[Theorem
1.7]{Pritsker*} has recently proved that no integer polynomial
$R(x)$
    can attain $\lh R(x)\rh_I=\tz(I)$,
    this value being achieved only by a normalized product of infinitely many
    polynomials. An immediate consequence of his
    result is the following.

\begin{prop} If an interval $I$ has an optimal monic integer Chebyshev polynomial then $\tm(I)>\tz(I)$.
\end{prop}

A fundamental question for both the monic and nonmonic integer transfinite
    diameter of an interval is whether its value can be computed exactly.
In \cite[Conjecture 5.1]{BorweinPinnerPritsker03}, Borwein {\em et al}
    make a conjecture for {\em Farey
    intervals} (intervals $\Farey$ where
    $b_1, b_2, c_1, c_2\in\BbZ$ and $b_2 c_1 - b_1 c_2=1$) concerning the exact value of their monic transfinite
     diameter.

\vspace{0.1 in} \noindent {\bf Conjecture BPP} (Farey Interval
Conjecture \cite[p. 82]{Borwein02}, \cite[Conjecture
5.1]{BorweinPinnerPritsker03}). 
\begin{em}
Suppose that $\Farey$ is a Farey interval, neither of whose endpoints is an integer.
Then
$$
\tm\left(\Farey\right)=\dfrac{1}{\min(c_1,c_2)}.
$$
\end{em}


Borwein {\em et al} verify their conjecture for all Farey intervals
    having  the denominators $c_1,c_2$  less than $22$.
In Section \ref{sec:Farey} we extend the verification to some
infinite families of Farey intervals
    (Theorems \ref{thm:farey} and \ref{thm:farey 2}).

 We next investigate what happens to $\tm([0, b])$ when $b$ is close
 to $\frac{1}{n}$.
For these intervals, some surprising things happen. Using the
polynomial $P(x)=x$, we know that $\tm([0,b])\le b<\frac{1}{n}$ if
$b<\frac{1}{n}$. In fact it appears likely that $\tm([0,b])$,
clearly a non-decreasing function of
    $b$, has a left discontinuity at $t=1/n\quad(n>1)$. On the
    other hand, we show in Theorem \ref{thm:interval 1/n} that $\tm$
    is locally constant on an interval of positive length  $\delta_n$ to the right of $\frac{1}{n}$.
Further, Theorem \ref{thm:interval 1/3} gives much larger values for
$\delta_n$
    for $n = 2, 3$ and $4$, as well as an upper bound for $\delta_2$.

In fact, more may be true.

\begin{conj}[Zero-endpoint Interval Conjecture]
\label{conj:[0,b]}
If $I = [0,b]$ is an interval with $b \leq 1$, then
    $\tm(I)=1/n$, where $n=\max\left(2,\left\lceil
    \frac{1}{b}\right\rceil\right)$
    is the smallest integer $n \geq 2$ for which $1/n \leq b$.
 \end{conj}
What little we know about $\tm([0,b])$ for $b>1$ is given in Theorem
\ref{thm:interval 1/3} (c), (d).

Both Conjecture BPP and Conjecture \ref{conj:[0,b]} are a consequence of the
    following conjecture.
\begin{conj}[Maximal obstruction implies $\tm(I)$ Conjecture]
\label{conj:maximal=tm} If an interval $I$ of length less than $4$
has a maximal obstruction $m$,
    then $\tm(I)=m$.
\end{conj}

We were at first tempted to conjecture here
    that $\tm(I)$, as well as equaling its maximal obstruction, is always attained by some monic integer polynomial.
However, the following counterexample eliminates this possibility in general.

\begin{cexample}
\label{thm:unattainable} The polynomial $7 x^3 + 4 x^2 - 2 x -1$ is
a maximal obstruction polynomial
    for the interval $I=[-0.684, 0.517]$.
However, there is no monic integer polynomial $P$ with $\lh P\rh_I$ equal to
    the maximal obstruction $7^{-1/3}$ for $I$.
\end{cexample}
This result is proved in Section \ref{sec:unattainable}.

Our next result proves the existence of maximal obstructions for
many intervals.
\begin{thm}
\label{thm:maximal} Every interval not containing an integer in its
interior has a maximal obstruction.
\end{thm}

 Based on
Conjecture \ref{conj:maximal=tm} and Theorem \ref{thm:maximal}
    we make the following conjecture.

\begin{conj}[Maximal Obstruction Conjecture]\label{conj:maximal}
Every interval of length less than $4$ has a maximal obstruction.
\end{conj}

We do not have much direct evidence for this conjecture. However,
our next conjecture, Conjecture \ref{conj:critical}, implies it.
    To describe this implication,
 we need the following notion, taken from Flammang, Rhin and Smyth
    \cite{FlammangRhinSmyth97}.
An irreducible polynomial $Q(x)=a_dx^d+\cdots +a_0\in \BbZ[x]$ with $a_d>0$, all
of whose
 roots lie in an interval $I$, and for which
 $a_d^{-1/d}$ is greater than the (nonmonic) transfinite diameter
 $\tz(I)$ is called a {\em critical polynomial} for $I$. Here we are
 interested only in nonmonic critical polynomials.

 It may be that every interval of length less than $4$ has infinitely many nonmonic critical polynomials -- see
  Proposition \ref{prop:crit} below. We make the following weaker conjecture.

\begin{conj}[Critical Polynomial Conjecture]\label{conj:critical}
    Every interval of length less than $4$ has at least one nonmonic
    critical polynomial.
\end{conj}

From Theorem \ref{thm:maximal} below, this conjecture is true for intervals not containing an integer. For intervals $I$ of length less than $4$ that do contain an integer (say $0$), then, since $\tz(I)<1$,  the polynomial $x$ is a critical polynomial for $I$. Thus `nonmonic' is an important word in this conjecture.

In Theorem \ref{thm:critical=maximal} we prove  that Conjecture \ref{conj:critical} implies Conjecture
    \ref{conj:maximal}.
More interestingly, we also prove in Corollary \ref{cor:critical=maximal} that Conjecture
\ref{conj:maximal=tm} and Conjecture \ref{conj:maximal} together imply
    Conjecture \ref{conj:critical}.

    We observe in passing the following conditional result for the integer transfinite diameter $\tz$.

\begin{prop}\label{prop:crit} Suppose that an interval $I$ has infinitely many critical polynomials
$Q_i(x)=a_{d_i,i}x^{d_i}+\dots+a_{0,i}$. Then
$$
\tz(I)=\inf_i a_{d_i,i}^{-\frac{1}{d_i}}.
$$
\end{prop}

This result is proved in Section \ref{sec:crit}. Montgomery \cite[p.182]{Montgomery94} conjectured this result
unconditionally for the interval $I=[0,1]$.


\section{Upper and Lower bounds for $\Lm(t)$ and $\Lp(t)$ for fixed $t$}
\label{sec:Upper/Lower Bounds}

The following lemma contains some simple properties, as well as
alternative definitions, of $\Lm$ and $\Lp$.

\begin{lem}
\label{lem:LL} We have
\begin{enumerate}
\item[(a)] $\Lm(t)\leq \Lp(t)$ for $t\ge 0$;
\item[(b)] $\Lm(t)=0$ for $0\le t \le \tfrac{1}{2}$;
\item[(c)]  $\Lp(t)\geq 2t$ for $0\leq t \le 1/2$;
 \item[(d)] $\Lm(t)=\sup_I\{d: \tm(I)\leq t \text{ for all } I \text{ with }
|I|=d\}$ for $t\geq \tfrac{1}{2}$;
\item[(e)] $\Lp(t)=\inf_I\{d: \tm(I)> t \text{ for all } I \text{ with }
 |I|=d\}$ for $t\ge 0$;
\item[(f)] $\Lp(t) = \Lm(t) = 4 t$ for $t \geq 1$.
\end{enumerate}
\end{lem}
\begin{proof}
First note that, by \cite[equation (1.11)]{BorweinPinnerPritsker03},
$\tm(I)=\tfrac{1}{2}$ for the zero-length interval
$\left[\tfrac{1}{2},\tfrac{1}{2}\right]$, from which (b) follows.

Part (c) follows from the fact that $\lh x\rh_{[-t,t]}=t$.

To prove (d), take $t\ge \tfrac{1}{2}$. Then the set
    \[S:=\{d: \tm(I)\leq t \text{ for all } I \text{ with } |I|=d\} \]
contains $0$ (by (b)), so
is nonempty. Put $s=\sup_d S$, and take $d\in S$. Since $I'\subset
I$ implies that $\tm(I')\le \tm(I)$ (\cite[Prop.
1.3]{BorweinPinnerPritsker03}), any $d'$ with $0\le d' < d$ also
lies in $S$, so that $S=[0,s)$ or $[0,s]$. Hence $\Lm(t)\ge s$. On
the other hand, for each $d>s$ there is an interval $I$ with $|I|=d$
and $\tm(I)>t$. Hence $\Lm(t)\le d$, giving $\Lm(t)=s$.

Now (a) follows straight from (b) and (d). The proof of (e),
similar to that of (d), is left as an exercise for the reader.

Finally, part (f) follows from the fact that for $|I| \geq 4$ we
have
    $\tm(I) = \tz(I) = \textrm{cap}(I) = \frac{|I|}{4}$
    (see for instance \cite{BorweinPinnerPritsker03}).
\end{proof}

 Next, we give a simple lemma,  needed for applying Proposition \ref{prop:L-bounds} below.

 \begin{lem}\label{L-simple} Suppose that $I_i=[a_i,b_i]\quad (i=1,\dots,n)$ are intervals with
 $a_1< a_2<\dots< a_n=a_1+1$,
   and put $M:=\max_{i=1}^{n-1}(b_{i+1}-a_i)$,
 $m:=\min_{i=1}^{n-1}(b_{i}-a_{i+1})$. Then
 \begin{enumerate}
 \item[(a)] Any interval of length at least $M$ contains an integer translate of some $I_i$.
 \item[(b)]  Any interval of length at most $m$ is contained in an integer translate of some $I_i$.
 \end{enumerate}
 \end{lem}

 \begin{proof}

  Given an interval $I$ of length $\ell$, we can, after translation by an integer, assume that $I=[a,b]$,
  where $a_j\le a< a_{j+1}$,  for some $j<n$.
 \begin{itemize}
 \item[(a)] Suppose that $\ell\ge M$. Then $b_{j+1}\le a_j+M\le a+\ell$, so that $[a_{j+1},b_{j+1}]\subset[a,a+\ell]$.
\item[(b)]  Suppose that $\ell\le m$. Then $b_{j}\ge a_{j+1}+m> a+\ell$, so that $[a,a+\ell]\subset[a_{j},b_{j}]$.
\end{itemize}
\end{proof}

The following proposition will be used to obtain explicit upper and
lower bounds for $\Lm(t)$ and $\Lp(t)$ for particular values of $t$.

\begin{prop}\label{prop:L-bounds} \noindent
\begin{enumerate}
\item[(a)] If $Q(x)=a_dx^d+\cdots+ a_0$, with integer coefficients and
$a_d>1$, has roots spanning an interval of length $\ell$, then for
any $t<a_d^{-1/d}$ we have
$$
\Lm(t)\le \ell.
$$
\item[(b)] Suppose that we have a finite set of polynomials
$Q_i(x)=a_{d_i,i}x^{d_i}+\cdots+a_{0,i}$ with all $a_{d_i,i}^{-1/d_i}>t$ with
    the property that every interval of length $\ell$ contains an integer
    translate of the roots of at least one of the polynomials $Q_i$. Then
$$
\Lp(t)\le \ell.
$$
\item[(c)] Suppose that we have a finite set of intervals $I_i$ such that
  for
each $I_i$ there is a monic integer polynomial $P_i$ with
$\lh P_i\rh_{I_i} \le t$. Suppose too that  every interval of length
$\ell$ is contained in an integer translate of some  $I_i$. Then
$$
\Lm(t)\ge \ell.
$$

\item[(d)]  If $\lh P\rh_I=t$ for some monic integer polynomial $P$ and
interval $I$ of length $\ell$, then
$$
\Lp(t)\ge \ell.
$$
\end{enumerate}
\end{prop}
\begin{proof} \noindent
\begin{enumerate}
\item[(a)] Given such a $Q(x), \ell$ and interval $I$ of length
$\ell$, and $t<a_d^{-1/d}$, then from Lemma BPP we have $\tm(I)\ge
a_d^{-1/d}>t$ so that, from the definition of $\Lm(t)$, we have $\Lm(t)\le
\ell$.

\item[(b)] Suppose that every interval $I$ of length $\ell$ contains
some integer translate of the set of roots of some $Q_i$. Then, by
Lemma BPP, $\tm(I)\ge a_{d_i,i}^{-1/d_i}>t$. Hence $\tm(I')>t$ for any
interval of length $|I'| \ge \ell$, and so $\Lp(t)\le\ell$.

\item[(c)] Here, for every interval $I$ of length $\ell$ with
$I+r\subset I_i$ say, (with $r \in \BbZ$), we have
\[ t>\lh P_i\rh_{I_i}\ge\lh P_i\rh_{I+r}=\lh P_i(x+r)\rh_I\ge\tm(I), \]
so that any $I'$ with $\tm(I')>t$ has $|I'|>\ell$. Hence $\Lm(t)\ge
\ell$.
\item[(d)] If $\lh P\rh_I=t$ and $|I|=\ell$ then $\tm(I)\le t$, so that $\Lp(t)\ge
\ell$.
\end{enumerate}
\end{proof}

\begin{table}[ht]
\begin{tabular}{|r|l|l|}
\hline
$i$ &Polynomials $Q_i$                        & Intervals $[a_i,b_i]$ \\
\hline
$1$ & $7x^3+7x^2-1 $ & $[-0.737, 0.328]$ \\
$2$ & $57x^6+81x^5+6x^4-32x^3-9x^2+3x+1 $ & $[-0.728, 0.494]$ \\
$3$ & $7x^3+4x^2-2x-1 $ & $[-0.684, 0.517]$ \\
$4$ & $59x^6+28x^5-43x^4-15x^3+11x^2+2x-1 $ & $[-0.669, 0.528]$ \\
$5$ & $3x^2-1 $ & $[-0.577, 0.577]$ \\
$6$ & $59x^6-28x^5-43x^4+15x^3+11x^2-2x-1 $ & $[-0.528, 0.669]$ \\
$7$ & $7x^3-4x^2-2x+1 $ & $[-0.517, 0.684]$ \\
$8$ & $57x^6-81x^5+6x^4+32x^3-9x^2-3x+1 $ & $[-0.494, 0.728]$ \\
$9$ & $7x^3-7x^2+1 $ & $[-0.328, 0.737]$ \\
$10$ & $63x^6-136x^5+72x^4+16x^3-17x^2+1 $ & $[-0.310, 1.115]$ \\
$11$ & $63x^6-146x^5+91x^4+7x^3-18x^2+x+1 $ & $[-0.285, 1.141]$ \\
$12$ & $58x^6-139x^5+90x^4+6x^3-18x^2+x+1 $ & $[-0.285, 1.178]$ \\
$13$ & $59x^6-147x^5+105x^4-3x^3-18x^2+2x+1 $ & $[-0.271, 1.184]$ \\
$14$ & $63x^6-159x^5+115x^4-4x^3-19x^2+2x+1 $ & $[-0.260, 1.197]$ \\
$15$ & $15x^4-29x^3+13x^2+x-1 $ & $[-0.244, 1.208]$ \\
$16$ & $57x^6-171x^5+153x^4-21x^3-21x^2+3x+1 $ & $[-0.228, 1.228]$ \\
$17$ & $15x^4-31x^3+16x^2-1 $ & $[-0.208, 1.244]$ \\
$18$ & $63x^6-219x^5+265x^4-126x^3+14x^2+5x-1 $ & $[-0.197, 1.260]$ \\
$19$ & $59x^6-207x^5+255x^4-127x^3+18x^2+4x-1 $ & $[-0.184, 1.271]$ \\
$20$ & $58x^6-209x^5+265x^4-136x^3+20x^2+4x-1 $ & $[-0.178, 1.285]$ \\
$21$ & $63x^6-232x^5+306x^4-171x^3+34x^2+2x-1 $ & $[-0.141, 1.285]$ \\
$22$ & $63x^6-242x^5+337x^4-204x^3+48x^2-1 $ & $[-0.115, 1.310]$ \\
\hline
\end{tabular}
\caption{Obstruction polynomials used for  Theorem
    \ref{thm:\Lm(1/2)} to prove that $\Lp(\tfrac{1}{2}) < 1.4715$.}
\label{tab:1.48}
\end{table}

\begin{table}[ht]
\begin{tabular}{|r|p{4in}|l|}
\hline $i$ & Polynomials $P_i$ & Intervals $I_i$ \\
\hline $1$ &$x^{1600} (x^3-4 x^2+1)^{36} (x^4+4 x^3-4 x^2-x+1)^{55}
\newline
 (x^8+236 x^7-96 x^6-167 x^5+64 x^4+39 x^3-14 x^2-3 x+1)^{39} \newline
 (x^8+372 x^7-196 x^6-249 x^5+129 x^4+55 x^3-28 x^2-4 x+2)^{20}
$  & $[-0.5142, 0.5613]$  \\&& \\

$2$ & $ x^{2121} (x^3-4 x^2+1)^{77} (x^4-10 x^3+5 x^2+2 x-1)^{84}
\newline (x^7-43 x^6-11 x^5+44 x^4+2 x^3-12 x^2+1)^{160}
$ & $[-0.4501, 0.5783]$ \\ && \\

$3$ & $ x^{12446} (x^2+x-1)^{199} (x^4-7 x^3+5 x^2+x-1)^{909}
\newline (x^6-53 x^5+46 x^4+10 x^3-14 x^2+1)^{640} $
& $[-0.4388, 0.5912]$ \\ && \\

$4$ & $ x^{312924} (x^4-7 x^3+5 x^2+x-1)^{45312}  \newline (x^4+8
x^3-8 x^2+1)^{217} (x^4+9 x^3-7 x^2-x+1)^{23800} $
& $[-0.4267, 0.6401]$ \\&& \\

$5$ & $x^{17556} (x^5+16 x^4-22 x^3 +5 x^2+3x-1)^{2256} \newline
(x^4+8 x^3-8 x^2+1)^{899} $
& $[-0.3797, 0.6847]$  \\&& \\

$6$ & $ x^{49329424964} (x-1)^{6557517120} (x^2+x-1)^{70328}\newline
(x^4+8 x^3-8 x^2+1)^{4916965515}  \newline (x^5-17 x^4+24 x^3-8
x^2-2 x+1)^{5952478752} \newline (x^5+16 x^4-22 x^3 +5
x^2+3x-1)^{541825536} $
& $[-0.3241, 0.7100]$ \\&& \\

$7$ & $x^{114080} (x-1)^{9324} (x^4+8 x^3-8 x^2+1)^{529} \newline
(x^4+9 x^3-9 x^2+1)^{2852} \newline (x^8+172 x^7-440 x^6+377 x^5-82
x^4-47 x^3+21 x^2+x-1)^{8184} \newline (x^8+214 x^7-531 x^6+440
x^5-90 x^4-54 x^3+23 x^2+x-1)^{6072} $
& $[-0.3064, 0.7344]$ \\&& \\

$8$ & $ x^{15200} (x-1)^{5192} (x^4+9 x^3-9 x^2+1)^{192} \newline
(x^8+172 x^7-440 x^6+377 x^5-82 x^4-47 x^3+21 x^2+x-1)^{1587}$  &
$[-0.2943, 0.7401]$ \\&& \\

$9$ & $ x^{3136} (x-1)^{1768} (x^6+3 x^5+6 x^4-18 x^3+9
x^2+x-1)^{32}
\newline (x^8+172 x^7-440 x^6+377 x^5-82 x^4-47 x^3+21 x^2+x-1)^{91}
$ & $[-0.2752, 0.7645]$  \\&& \\

$10$ & $x^{146704} (x-1)^{85868} (x^2+x-1)^{6369}  \newline
(x^6+3x^5+6x^4-18x^3+9x^2+x-1)^{1768}$ & $[-0.2622, 1.1030]$  \\
\hline
\end{tabular}
\caption{Optimal monic integer Chebyshev polynomials used for
Theorem
    \ref{thm:\Lm(1/2)} to prove that $\Lm(\tfrac{1}{2})\ge 1.008848$.}
\label{tab:interval 1}
\end{table}

\clearpage

\begin{proof}[Proof of Theorem \ref{thm:\Lm(1/2)}]

Applying Proposition \ref{prop:L-bounds}(a) with $Q(x)=7 x^3-7
x^2+1$, we have
$$\Lm(\tfrac{1}{2})\le\ell=1.064961507.$$
Here, a more precise value could be determined by calculating the
    span of the roots of $Q(x)$ to a higher precision.

We apply Proposition \ref{prop:L-bounds}(b) and Lemma
\ref{L-simple}(a) using the polynomials $Q_i$ of
    Table \ref{tab:1.48}, with the intervals $[a_i,b_i]$ containing their roots.
(Here, the endpoints listed in Table \ref{tab:1.48} are
approximations
    of the minimal and maximal root of the obstruction polynomial in question.
A higher precision was used for the computation of the upper bound of
    $\Lp\left(\tfrac{1}{2}\right) < 1.4715$.)
We put $Q_{23}(x)=Q_1(x-1)$, whose  roots are contained
    in $[a_{23},b_{23}]:=[a_1+1,b_1+1]$, and apply the Proposition to
    the $23$ polynomials $Q_1,\cdots,Q_{23}$.
Each has $a_d^{-1/d}>\tfrac{1}{2}$. Then because $\max_{i=1}^{22}
(b_{i+1}-a_i)=b_{16}-a_{15}=1.4715$, any interval $I$ of length $|I|
>1.4715$ must, by Lemma \ref{L-simple}(a), contain  some integer
translate of some interval $[a_i,b_i]$, and so all the
 roots of  the corresponding polynomial $Q_i$. Hence $\Lp(\tfrac{1}{2})<1.4715$.

We apply Proposition \ref{prop:L-bounds}(c) by starting with the
$10$ intervals $I_i\quad(i=1,\cdots,10)$ in Table \ref{tab:interval
1}, and putting
 $I_i=1-I_{21-i}$ and $P_i(x)=P_{21-i}(1-x)$ for $i=11,\cdots,20$, with $I_{21}=1+I_1$ and $P_{21}(x)=P_1(x-1)$.
(Here again, the endpoints listed in Table \ref{tab:interval 1} are
    approximations only.  To find a more accurate values, we would solve
    for the roots of $P(x) = \pm \left(\frac{1}{2}\right)^{\deg P}$.
Higher precision values were used to compute the lower bound
    $\Lm\left(\tfrac{1}{2}\right) > 1.008848$.)
Each polynomial $P_i$ listed has a critical point at $\tfrac{1}{2}$
    (and also at $-\tfrac{1}{2}$ in the case of the last polynomial),
    with $P_i\left(\tfrac{1}{2}\right)=\pm
    \left(\tfrac{1}{2}\right)^{\deg{P_i}}$.
The value of $|P_i(x)|$ at all other critical points, as well as at
    the interval endpoints, is strictly less than
    $\left(\tfrac{1}{2}\right)^{\deg{P_i}}$.
This shows in each case that  $\lh P_i\rh_{I_i}= \tfrac{1}{2}$.
Then all $21$ intervals $I_i$ have $\tm(I_i)=\tfrac{1}{2}$
    and, writing $I_i=[a_i,b_i]\quad(i=1,\cdots,21)$ we have
 \begin{equation}\label{E:20}
 \min_{i=1}^{20}\left(b_{i}-a_{i+1}\right)=b_5-a_6>1.008848.
 \end{equation}
   From this it follows by Lemma \ref{L-simple}(b) that every interval $I$ of length less than $1.008848$ is a
 subinterval of an integer translate of some $I_i$, so that $\tm(I)\le\lh P_i\rh_{I}\le \lh P_i\rh_{I_i}= \tfrac{1}{2}$.
 This proves part (a) of the Theorem.

Part (b) of the Theorem follows on applying Proposition
\ref{prop:L-bounds} (d) with $P(x)=x^2-x$.  We then have, for
$I=\left[\frac{1-\sqrt{2}}{2}, \frac{1+\sqrt{2}}{2}\right]$, that
 $\lh P\rh_I=\tfrac{1}{2}$, so that $\Lp(t)\ge
|I|=\sqrt{2}$.

\end{proof}

\begin{table}[ht]
\begin{tabular}{|l|l|l|l|}
\hline
$i$ & Polynomial $Q_i$                          & $t_i$         &  $\ell_i^-$         \\
\hline
1 & $7 x^3+7 x^2-1                         $
  & $\frac{1}{\sqrt[3]{7}} \approx 0.522$ & 1.064961507  \\
2 & $3 x^2-1                               $
  & $\frac{1}{\sqrt[2]{3}} \approx 0.577$ & 1.154700538  \\
3 & $5 x^3+3 x^2-2 x-1                     $
  & $\frac{1}{\sqrt[3]{5}} \approx 0.584$ & 1.390656045  \\
4 & $2 x^2-1                               $
  & $\frac{1}{\sqrt[2]{2}} \approx 0.707$ & 1.414213562  \\
5 & $3 x^4-2 x^3-4 x^2+x+1                 $
  & $\frac{1}{\sqrt[4]{3}} \approx 0.759$ & 2.173182852  \\
6 & $2 x^3-4 x^2+1                         $
  & $\frac{1}{\sqrt[3]{2}} \approx 0.793$ & 2.306243643  \\
7 & $2 x^4-8 x^3+8 x^2-1                   $
  & $\frac{1}{\sqrt[4]{2}} \approx 0.840$ & 2.613125930  \\
8 & $2 x^5-15 x^4+39 x^3-40 x^2+12 x+1     $
  & $\frac{1}{\sqrt[5]{2}} \approx 0.870$ & 2.982466529  \\
9 & $2 x^6-12 x^5+22 x^4-8 x^3-10 x^2+4 x+1$
  & $\frac{1}{\sqrt[6]{2}} \approx 0.890$ & 3.131521012  \\
\hline
\end{tabular}
\caption{Upper bounds for $\Lm(t)$. Here $\Lm(t)<\ell_i^-$ for $t<t_i$,
where $\ell_i^-$ is the span of the roots of the $i$th polynomial (see
Theorem \ref{thm:\Lm(t)}).} \label{tab:\Lm(t)}
\end{table}

\begin{table}[ht]
\begin{tabular}{|r|c|c|| c|c|c|}
\hline
$i$ & $t_i$ & $\ell_i^+$ & $i$ & $t_i$ & $\ell_i^+$ \\
\hline
1 & $ \frac{1}{\sqrt[6]{63}} \approx .501 $ & 1.47149 &
31 & $ \frac{1}{\sqrt[5]{15}} \approx .582 $ & 1.71707 \\
2 & $ \frac{1}{\sqrt[6]{60}} \approx .505 $ & 1.47887 &
32 & $ \frac{1}{\sqrt[3]{5}} \approx .585 $ & 1.72578 \\
3 & $ \frac{1}{\sqrt[5]{30}} \approx .506 $ & 1.48183 &
33 & $ \frac{1}{\sqrt[6]{24}} \approx .589 $ & 1.78511 \\
4 & $ \frac{1}{\sqrt[6]{59}} \approx .507 $ & 1.48424 &
34 & $ \frac{1}{\sqrt[5]{14}} \approx .590 $ & 1.79006 \\
5 & $ \frac{1}{\sqrt[4]{15}} \approx .508 $ & 1.48823 &
35 & $ \frac{1}{\sqrt[6]{23}} \approx .593 $ & 1.80103 \\
6 & $ \frac{1}{\sqrt[6]{58}} \approx .508 $ & 1.49541 &
36 & $ \frac{1}{\sqrt[4]{8}} \approx .595 $ & 1.80333 \\
7 & $ \frac{1}{\sqrt[6]{57}} \approx .510 $ & 1.49802 &
37 & $ \frac{1}{\sqrt[6]{22}} \approx .597 $ & 1.80514 \\
8 & $ \frac{1}{\sqrt[6]{56}} \approx .511 $ & 1.50442 &
38 & $ \frac{1}{\sqrt[6]{19}} \approx .612 $ & 1.82308 \\
9 & $ \frac{1}{\sqrt[4]{14}} \approx .517 $ & 1.50918 &
39 & $ \frac{1}{\sqrt[4]{7}} \approx .615 $ & 1.82808 \\
10 & $ \frac{1}{\sqrt[6]{51}} \approx .519 $ & 1.51232 &
40 & $ \frac{1}{\sqrt[6]{18}} \approx .618 $ & 1.85414 \\
11 & $ \frac{1}{\sqrt[3]{7}} \approx .523 $ & 1.51409 &
41 & $ \frac{1}{\sqrt[5]{11}} \approx .619 $ & 1.86446 \\
12 & $ \frac{1}{\sqrt[6]{48}} \approx .525 $ & 1.54721 &
42 & $ \frac{1}{\sqrt[6]{17}} \approx .624 $ & 1.86909 \\
13 & $ \frac{1}{\sqrt[5]{25}} \approx .525 $ & 1.54825 &
43 & $ \frac{1}{\sqrt[3]{4}} \approx .630 $ & 1.87806 \\
14 & $ \frac{1}{\sqrt[4]{13}} \approx .527 $ & 1.55329 &
44 & $ \frac{1}{\sqrt[6]{15}} \approx .637 $ & 1.92375 \\
15 & $ \frac{1}{\sqrt[6]{46}} \approx .528 $ & 1.56522 &
45 & $ \frac{1}{\sqrt[5]{9}} \approx .644 $ & 1.92862 \\
16 & $ \frac{1}{\sqrt[6]{45}} \approx .530 $ & 1.57021 &
46 & $ \frac{1}{\sqrt[4]{5}} \approx .669 $ & 1.95815 \\
17 & $ \frac{1}{\sqrt[5]{23}} \approx .534 $ & 1.57066 &
47 & $ \frac{1}{\sqrt[6]{11}} \approx .671 $ & 2.03528 \\
18 & $ \frac{1}{\sqrt[4]{12}} \approx .537 $ & 1.57390 &
48 & $ \frac{1}{\sqrt[3]{3}} \approx .693 $ & 2.05072 \\
19 & $ \frac{1}{\sqrt[5]{21}} \approx .544 $ & 1.58148 &
49 & $ \frac{1}{\sqrt{2}} \approx .707 $ & 2.07313 \\
20 & $ \frac{1}{\sqrt[4]{11}} \approx .549 $ & 1.59285 &
50 & $ \frac{1}{\sqrt[6]{7}} \approx .723 $ & 2.46521 \\
21 & $ \frac{1}{\sqrt[5]{20}} \approx .549 $ & 1.60583 &
51 & $ \frac{1}{\sqrt[5]{5}} \approx .725 $ & 2.49418 \\
22 & $ \frac{1}{\sqrt[6]{36}} \approx .550 $ & 1.62320 &
52 & $ \frac{1}{\sqrt[6]{6}} \approx .742 $ & 2.55291 \\
23 & $ \frac{1}{\sqrt[6]{34}} \approx .556 $ & 1.63662 &
53 & $ \frac{1}{\sqrt[5]{4}} \approx .758 $ & 2.58796 \\
24 & $ \frac{1}{\sqrt[6]{33}} \approx .558 $ & 1.64392 &
54 & $ \frac{1}{\sqrt[4]{3}} \approx .760 $ & 2.60202 \\
25 & $ \frac{1}{\sqrt[5]{18}} \approx .561 $ & 1.65596 &
55 & $ \frac{1}{\sqrt[6]{5}} \approx .765 $ & 2.61238 \\
26 & $ \frac{1}{\sqrt[6]{32}} \approx .561 $ & 1.65815 &
56 & $ \frac{1}{\sqrt[6]{4}} \approx .794 $ & 2.70928 \\
27 & $ \frac{1}{\sqrt[4]{10}} \approx .562 $ & 1.66032 &
57 & $ \frac{1}{\sqrt[5]{3}} \approx .803 $ & 2.89569 \\
28 & $ \frac{1}{\sqrt[6]{31}} \approx .564 $ & 1.66308 &
58 & $ \frac{1}{\sqrt[6]{3}} \approx .833 $ & 2.97756 \\
29 & $ \frac{1}{\sqrt[5]{16}} \approx .574 $ & 1.67218 &
59 & $ \frac{1}{\sqrt[4]{2}} \approx .841 $ & 2.98928 \\
30 & $ \frac{1}{\sqrt{3}} \approx .577 $ & 1.68244 &
60 & $ \frac{1}{\sqrt[5]{2}} \approx .871 $ & 3.23520 \\
\hline
\end{tabular}
\caption{Upper bounds for $\Lp(t)$. \newline Here $\Lp(t)<\ell_i^+$ for
$t<t_i$, where $\ell_i^+$ is the span of the roots of the $i$th
polynomial (see Theorem \ref{thm:\Lm(t)}).} \label{tab:\Lp(t)}
\end{table}


\section{General bounds for $\Lm(t)$ and $\Lp(t)$}

In this section we find upper and lower bounds for $\Lm(t)$ and $\Lp(t)$,
 valid for $t$ from $0.5$ to close to $0.9$. Our first result gives the upper bounds.

\begin{thm} \label{thm:\Lm(t)}  \noindent
\begin{itemize}
\item[(a)] For all $t_i$ and
$\ell_i^-$ in Table \ref{tab:\Lm(t)} and for
    all $t<t_i$ we have $\Lm(t) < \ell_{i}^-$.
 \item[(b)]For all $t_i$ and $\ell_i^+$ in Table
\ref{tab:\Lp(t)} and for
    all $t<t_i$ we have $\Lp(t) < \ell_{i}^+$.
\end{itemize}
\end{thm}

The Theorem is proved by applying Proposition \ref{prop:L-bounds} (a) and
    (b) for a range of values in $[0.5,1]$.
Here again, the diameter given in Table \ref{tab:\Lm(t)} can be
    computed more exactly by considering the difference
    between the maximal and minimal roots of the obstruction polynomial.
For Table \ref{tab:\Lp(t)}, a calculation similar to that done for Table
    \ref{tab:1.48} was done for each $t_i$.
The  rounding  procedure was that used for Table \ref{tab:1.48}.
Then the monotonicity of $\Lm(t)$ and $\Lp(t)$ gives the result for
    all $t$ in this range.

For the lower bounds, we first define the normalized polynomial $P_\alpha$
 \begin{equation}\label{eqn:normal}
 P_\alpha(x)
= (x (1-x))^{\tfrac{1-\alpha}{2}}
    (x^2-x-1)^{\tfrac{\alpha}{2}},
\end{equation}
     of degree $1$, and let $\al^*\approx 0.4358$ be the root in $(0,1)$ of the
equation \begin{equation}\label{eqn:root}
4\al^\al(1-\al)^{1-\al}=5^\al.
\end{equation}

The following result gives the lower bounds.

\begin{prop}\label{P-lowerB}
For $0\le\al\le \frac{\ln 4}{\ln 5}$ we have
\begin{itemize}
\item[(a)]
 $\Lp\left(\frac{5^{\alpha/2}}{2}\right)\ge \ell_\al$, where $\ell_\al$ is the
root of $P_\al\left(\tfrac{1}{2}+\ell_\al/2\right)=\frac{5^{\alpha/2}}{2}$ in
\[\begin{cases} (\sqrt{5},\infty) &\text{ if } \al>\al^*;\\
               (1,\sqrt{5}) &\text{ if } \al\le\al^*.
\end{cases}
\]
\item[(b)] $\Lm\left(\frac{5^{\alpha/2}}{2}\right)\ge \max(\ell_\al-1,1.008848)$.
\end{itemize}
\end{prop}

\begin{figure}[ht]
\[\psfig{file=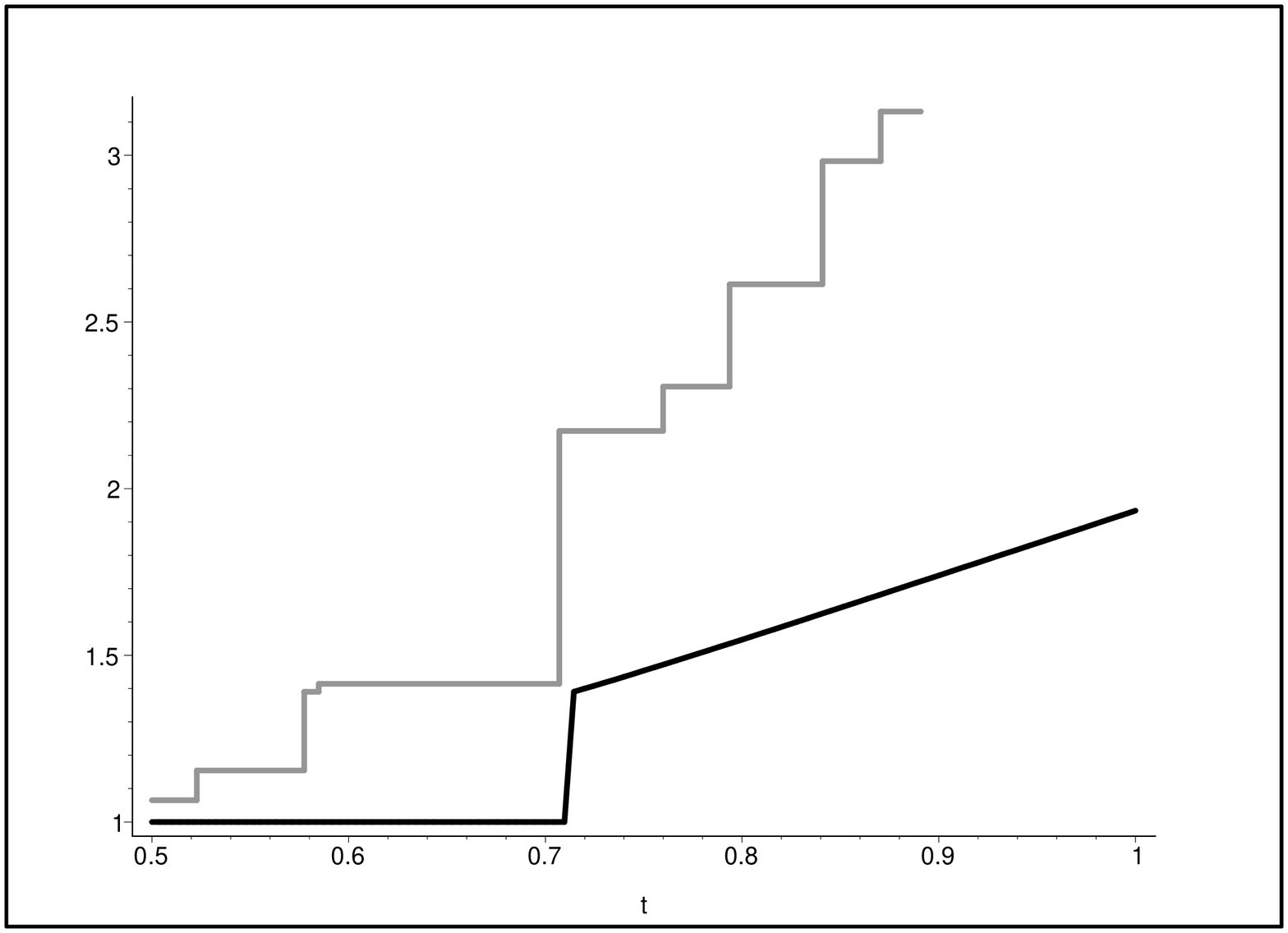,height= 300pt,width= 300pt,angle=270}\]
\caption{Upper and lower bounds for $\Lm(t)$ (Theorem
\ref{thm:\Lm(t)} and
    Proposition \ref{P-lowerB}). \newline
         grey line -- upper bound; \newline
         black line -- lower bound.}
\label{fig:\Lm(t)}
\end{figure}

\begin{figure}[ht]
\[\psfig{file=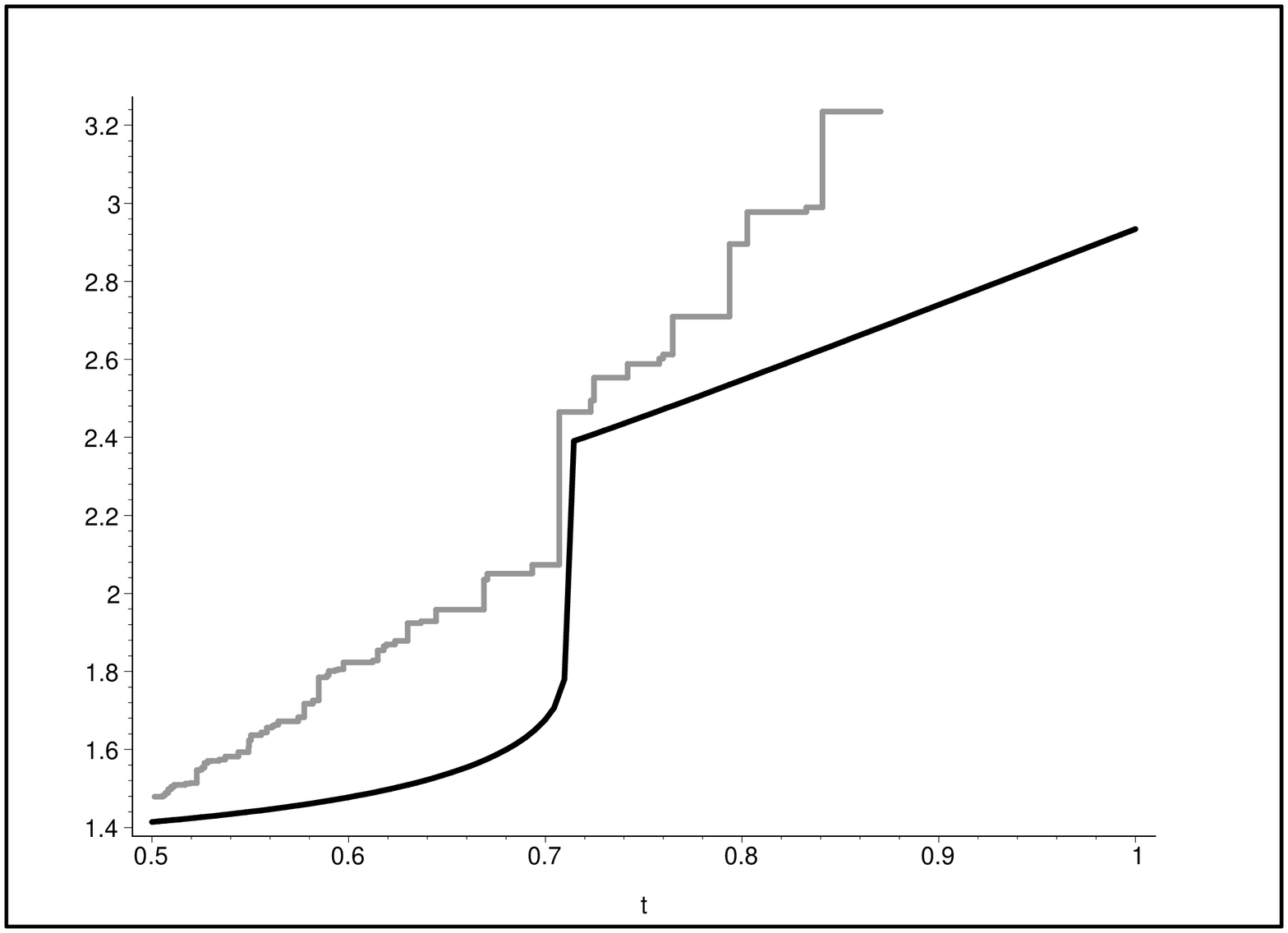,height= 300pt,width= 300pt,angle=270}\]
\caption{Upper and lower bounds for $\Lp(t)$ (Theorem
\ref{thm:\Lm(t)} and Proposition \ref{P-lowerB}). \newline
         grey line -- upper bound;  \newline
         black line -- lower bound.}
\label{fig:\Lp(t)}
\end{figure}


   For the proof, we need the following simple observation.

\begin{lem}
\label{lem:Ll} If $\Lp(t) \geq \ell+1$ then $\Lm(t) \geq \ell$.
\end{lem}

This follows straight from the fact that, given an interval $I$ of
length $\ell+1$, every interval of length $\ell$ has an integer translate
that is a subinterval of $I$.

\begin{proof}[Proof of Proposition \ref{P-lowerB}]
It should first be pointed out that this proposition is in fact true
   for all $\alpha$, and not just those in the range specified.
That being said, for $\alpha > \frac{\ln 4}{\ln 5}$ we would have
    $\frac{5^{\alpha/2}}{2} > 1$, in which case we could appeal to
    Lemma \ref{lem:LL} (f) for the exact answer.
\begin{itemize}
\item[(a)]
We will proceed to analyze $\lh P_\alpha(x)\rh_{I_\ell}$,
    picking $\alpha$ and
    $\ell$ such that, at the endpoints of the interval $I_\ell$, $|P\al(x)|$ equals the largest local maximum
     of $|P\al(x)|$ in the interior of $I_\ell$.
(Notice that $P_\al$ is already normalized, so $\lh P_\al \rh_{I_\ell} = ||P_\al||_{I_\ell}||$.)
(See Figures \ref{fig:alpha<alpha*} and \ref{fig:alpha>alpha*}.)

\begin{figure}[ht]
\[\psfig{file=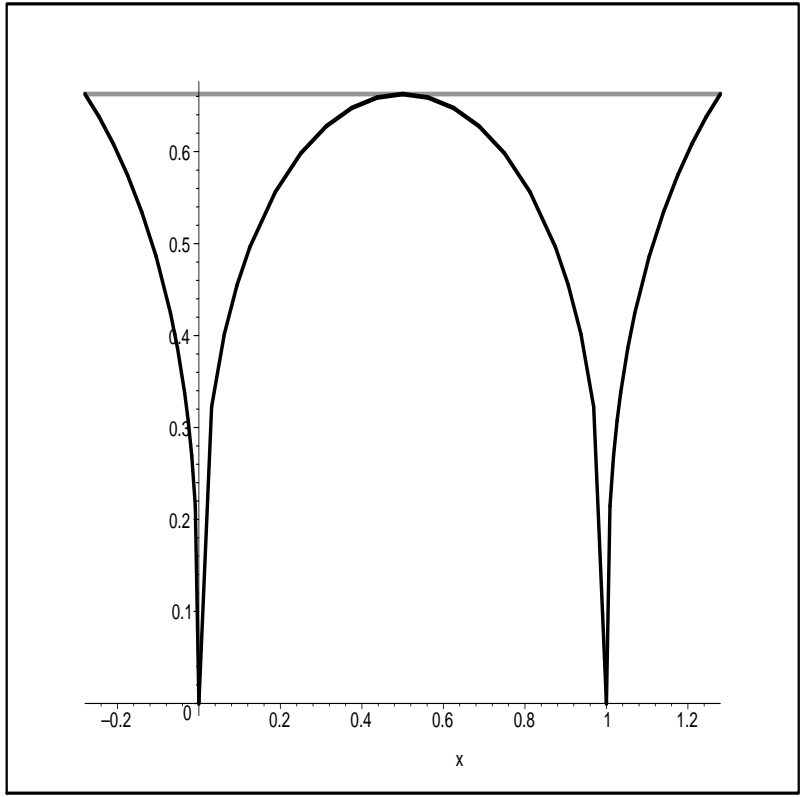,height= 300pt,width= 300pt,angle=270}\]
\caption{The normalized polynomial $P_\al(x)$ (see (\ref{eqn:normal}))
with $\alpha = 0.35 < \alpha^*$,
         $\ell_\alpha \approx 1.559$ and
         $\tm(I_{\ell_\al})\le \lh P_\alpha(x)\rh_{I_{\ell_\al}} \approx 0.663$.}
\label{fig:alpha<alpha*}
\end{figure}

\begin{figure}[ht]
\[\psfig{file=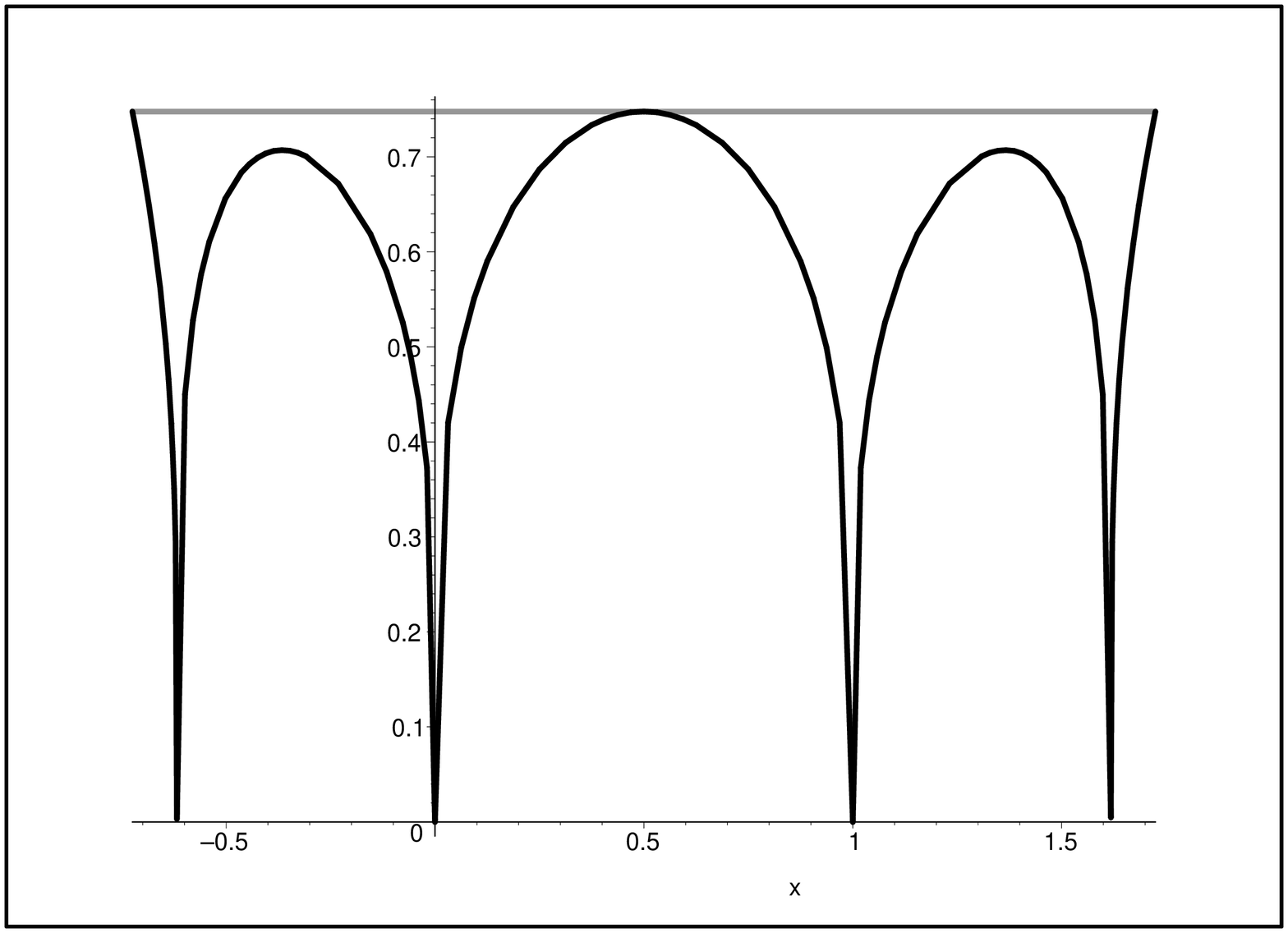,height= 300pt,width= 300pt,angle=270}\]
\caption{The normalized polynomial $P_\al(x)$ (see (\ref{eqn:normal}))
with $\alpha = 0.5 > \alpha^*$,
         $\ell_\alpha \approx 2.449$ and
         $\tm(I_{\ell_\al})\le \lh P_\alpha(x)\rh_{I_{\ell_\al}} \approx 0.748$.}
\label{fig:alpha>alpha*}
\end{figure}


Notice first that
\[
    \left|P_\alpha\left(\tfrac{1}{2}+x/2\right)\right|
    =\frac{5^{\tfrac{\alpha}{2}}}{2}\left|1-x^2\right|^{\tfrac{1-\alpha}{2}}
    \left|1-\frac{x^2}{5}\right|^{\tfrac{\alpha}{2}},
 \]
which has a local maximum of  $\frac{5^{\alpha/2}}{2}$ at $x=0$, and
a local maximum of
$m_\alpha=|1-\alpha|^{(-\alpha)/2}|\alpha|^{\alpha/2}$ at
$x^2=5-4\alpha$ . Now the equation $m_\al=\frac{5^{\alpha/2}}{2}$
has a root $\al$ defined by (\ref{eqn:root}), with
$m_\al>\frac{5^{\alpha/2}}{2}$ for $\al<\al^*$ and
$m_\al<\frac{5^{\alpha/2}}{2}$  for $\al>\al^*$. Hence if $\al\ge
\al^*$ then $|P_\alpha(\tfrac{1}{2}+x/2)|\le\frac{5^{\alpha/2}}{2}$
for $x\le \sqrt{5}$, so that
$\lh P_\al\rh_{I_\al}=\frac{5^{\alpha/2}}{2}$, where
$I_\al=\left[\tfrac{1}{2}-\ell_\al/2,\tfrac{1}{2}+\ell_\al/2\right]$ with
$\ell_\al$ the root $\ell_\al>\sqrt{5}$ of
$P_\alpha\left(\tfrac{1}{2}+\ell_\al/2\right) =\frac{5^{\alpha/2}}{2}$.
However, if $\al<\al^*$ then we have the same result, but only for
$\ell_\al$ the
 root in $(1,\sqrt{5})$ of $P_\alpha(\tfrac{1}{2}+\ell_\al/2)=\frac{5^{\alpha/2}}{2}$. This
  gives the lower bound $\Lp\left(\frac{5^{\alpha/2}}{2}\right)\ge \ell_\al$,
   but with a
  left discontinuity in $\ell_\al$ (as a function of $\al$) at $\al=\al^*$.
    A plot of this lower bound, along with the upper bounds from
    Theorem \ref{thm:\Lm(t)} and Table \ref{tab:\Lp(t)}, is given
    in Figure \ref{fig:\Lp(t)}.

\item[(b)] We know that $\Lm$ is a non-decreasing function,
    and that $\Lm(\tfrac{1}{2}) \geq 1.008848$. Combining these facts with Lemma \ref{lem:Ll} we get that
    $\Lm\left(\frac{5^{\alpha/2}}{2}\right)
    \geq \max(\ell_\alpha - 1, 1.008848)$.
This is displayed numerically, along with the upper bounds from
    Theorem \ref{thm:\Lm(t)} and Table \ref{tab:\Lm(t)},
    in Figure \ref{fig:\Lm(t)}.

\end{itemize}
\end{proof}


\section{Intervals of length 1: Proof of Theorem \ref{thm:interval 1}}
\label{sec:interval 1}

\begin{proof}[Proof of Theorem \ref{thm:interval 1}] From Theorem
 \ref{thm:\Lm(1/2)} (a) we know that
  every interval $I$ of length $\ell \le 1.008848$ has $\tm(I)\le
 \frac{1}{2}$.
Now since  every interval of length $\ell \ge 1$ has some integer
translate that contains $\tfrac{1}{2}$, we have
 $$
 \tfrac{1}{2}=\tm\left(\left\{\tfrac{1}{2}\right\}\right)\le \tm(I)
 $$ for all such intervals, so that $\tm(I)=\tfrac{1}{2}$ for all $I$ with
 $1\le |I|\le 1.008848$.

 If $b>1.064961507$ then again from Theorem \ref{thm:\Lm(1/2)} (a), with the
    polynomial $Q(x) = 7 x^3 + 7 x^2 - 1$, there is an interval $I$ of
    length $b$ with $\tm(I)>\frac{1}{2}$.

To complete the proof, note that for $b<1$
\begin{align*}
\tm([-b/2,b/2])&=\sqrt{\tm([0,b^2/4])},\\
\intertext{on applying \cite[Prop 1.4 with the polynomial
$x^2$]{BorweinPinnerPritsker03},  and then }
\sqrt{\tm([0,b^2/4])} & \le \sqrt{b^2/4}<\tfrac{1}{2},
\end{align*}
using the polynomial $P(x)=x$ on $[0,b^2/4]$.
\end{proof}

\section{Computational methods}\label{sec:comp}

\subsection{Finding optimal monic integer Chebyshev polynomials {\em P}} We now describe  how
 the polynomials of Table \ref{tab:interval 1}
were found.  These are optimal monic integer polynomials $P$ having $\lh P\rh_I=\tfrac{1}{2}$ on various intervals of
length just greater than $1$. For these intervals, the maximal obstruction polynomial is $Q(x)=2x-1$, and the
maximal obstruction is $m=\tfrac{1}{2}$. The method applies more generally, however, to any interval $I$ having
a maximal obstruction polynomial $Q$, so we shall describe the method for this more general situation. We suppose
 that the maximal obstruction is $m=a_d^{-1/d}$, where
$Q(x)=a_dx^d+\dots+a_0$, so that we seek a monic integer polynomial
$P$ with $\lh P\rh_I=m$.

 Firstly, potential factors of $P$ of small degree $k$ were
identified using LLL
    \cite{BorweinPinnerPritsker03, Hare02a, LenstraLenstraLovasz82}.
The basis used was $[1, x, \cdots, x^k]$, with the inner product
\[
    \langle R_1, R_2 \rangle = \int_I R_1(x) R_2(x)\ dx + b_k c_k.
\]
    Here $R_1(x) = b_k x^k + \cdots + b_0$ and
     $R_2(x) = c_k x^k + \cdots + c_0$.
The $b_k c_k$ component of the inner product was inserted to
    discourage nonmonic
    polynomials from appearing in the basis returned by LLL.
Now, at least one element in the basis will contain
    an $x^k$ term and,
    because of the $b_k c_k$ penalty, such an element is almost always
monic. (In fact always in the examples we computed.) So we obtained
a monic polynomial of degree $d$ with small $L_2$ norm, which
usually also had a small supremum norm. These monic polynomials with
small $L_2$ norm are not necessarily
    irreducible.
At this point we examined each of their irreducible factors $f_i$, again monic polynomials,
  and applied Lemma \ref{lem:resultant-value}(a) below to eliminate
    some of them.
We then used the method of Borwein and Erd\' elyi
\cite{BorweinErdelyi96} to search for exponents $\al_i \in \BbN$
such that $P^{1/\deg P}:=\prod_i f_i^{\al_i/\deg f_i}$
    has the desired property $\lh P\rh_I=m$.
To do this, we needed to minimize $t$ subject to the constraint
    \[ \sum_i \frac{\alpha_i}{\deg f_i} \log(|f_i(x)|) \leq t,\] for all
    $x \in I$ with $\sum_i \alpha_i = 1$,
    $0 \leq \alpha_i$.
Some additional constraints on the $\alpha_i$ that we made use of
are given by Lemma \ref{lem:resultant-value} (b), (c). The main difference
between our application and the original one is that
    here the  polynomials $f_i$ are all monic.
By choosing a large number of points $x \in I$ to discretize the
problem, we get a system of linear equations, on which
    the Simplex method can be used to get a good estimate of $\min(t)$
    \cite{BorweinErdelyi96, HabsiegerSalvy97, Schrijver86}.
In practice, with a high enough precision and a large enough number
of    sample points, we obtain $\min(t)=m$ exactly, and the
corresponding $\al_i$ then give the required $P$. We then check that
$P$
    is indeed an optimal monic integer Chebyshev polynomial for $I$ by checking algebraically that
    $|P|^{1/\deg P}= m $ at all roots of the maximal obstruction polynomial $Q$, and
    furthermore   that all other local
    maxima of $|P|$ in this interval are strictly smaller than $m$.

The following lemma, used to help construct these polynomials $P$,
specifies extra properties that their  factors $f_i$ and normalized
exponents $\al_i$ must have.

\begin{lem}
\label{lem:resultant-value}  Let $I$ be an interval that has  a
maximal obstruction polynomial $Q(x) = a_d x^d + \cdots +a_0$.
Suppose further that $P(x)$ attains the maximal obstruction, and
that $P(x)^{1/\deg P}
    =\prod_if_i^{\alpha_i/\deg f_i}$, with $\sum_i\alpha_i=1$.
Then
\begin{enumerate}
\item[(a)]
 The resultant $\Res(f_i, Q)$ is equal to $ \pm 1$ for every factor $f_i\in\BbZ[x]$
of $P$.

\item[(b)] For every root $\beta$ of $Q$ we have
$$
\sum_i\frac{\al_i}{\deg f_i} \times \frac{f'_i(\beta)}{f_i(\beta)}=0.
$$

\item[(c)] Fix a root $\beta\in\mathbb R$ of $Q$, and put $\hat f_i=|f_i(\beta)|^{1/\deg f_i}\in\mathbb R$.
     Let $\F$ be the
    multiplicative subgroup of $\mathbb R_{>0}$ generated
    by $a_d$ and the $\hat f_i$
with $b_1=a_d$ and $b_2, \cdots, b_k$ an independent generating set
for $\F$, with say ${\hat f_i}^{1/\deg f_i} = \prod_j b_j^{ c_{j,i}
}$ for some integers $c_{j,i}$. Then
\[ \sum_i c_{j,i} \alpha_i =\begin{cases}
 -1/d &\text{ if } j=1;\\
  0 &\text{ if } j>1.
\end{cases}
\]
\end{enumerate}
\end{lem}

\begin{proof}

We have  $\prod_i P(\beta_i) = \pm 1/a_d^{\deg P}$, where the
product is taken over the    roots $\beta_i$ of $Q$, so that, from
(\ref{eqn:resultant}), $\Res(P, Q) = \pm 1$. Then (a) follows from
the fact that the resultant of a product with $Q$ is the product of the
    resultants with $Q$.

The second part follows from the fact that all  the roots $\beta$ of
$Q$ must be critical points of $P(x)$. Further, since $P(x)$ attains
the maximal obstruction, we  have from Lemma BPP that for all such
$\beta$ we have $|P(\beta)|^{1/\deg P} =  a_d^{-1/d}$, giving the
third part.
\end{proof}

Note that Lemma \ref{lem:resultant-value}
     simplifies considerably when the maximal obstruction polynomial
     is linear, say $a_1 x - a_0$.
Then it says that
    $f_i\left(\frac{a_0}{a_1}\right) = \pm a_1^{-\deg f_i}$ and
    with  $P'\left(\frac{a_0}{a_1}\right) = 0$.

    The independent generating set $b_1,\cdots,b_k$ for $\F$ was found using
    the integer relation-finding program PSLQ, which we used to
    search for linear integer relations between $\log a_d$ and the
    $\log\hat f_i$.

 As we have seen, the method for finding an optimal monic integer Chebyshev polynomial $P$ depends on first finding the
 (in practice there was only one) maximal obstruction polynomial for the interval. We now describe how to do this.

\subsection{Finding obstruction polynomials {\em Q}}

The obstruction polynomial $7 x^3-7 x^2+1$, as well as those listed in Table
\ref{tab:1.48}
    and \ref{tab:\Lm(t)},
    were found using the technique of Robinson \cite{Robinson64}
     (see also \cite{McKeeSmyth04, Smyth84}).
In this method, the aim is to search for all degree $d$ polynomials
$Q(x)=a_dx^d+\cdots +a_0$
    having all their roots in an interval $I_0$, for fixed degree,
    and fixed lead coefficient, $a_d$, with $a_d \leq 2^d$.
We describe below how $I_0$ is chosen. Robinson's method uses the
fact that for $k=1, 2, \cdots, d-1$ the span of
    the roots of the $k$th derivative of $Q$ is contained in the span of
    the roots of the $(k-1)$th derivative of $Q$.
In particular, these derivatives  have all their roots in $I_0$.

Starting with the $(d-1)$st derivative of $Q$, we get a range of
possible valid values
    for $a_{d-1}$.
Consider then the $(d-2)$nd derivative
    to find valid ranges for $a_{d-2}$.
Continuing in this fashion, we obtain a list of polynomials, each one having
all its roots
    in $I_0$.
We now sieve this list, first by eliminating all polynomials that are reducible, or have
integer content
    greater than $1$.
Having obtained a list of irreducible polynomials, we can then prune
it
    further, as follows.
If $Q(x)$ and $R(x)$ are both irreducible polynomials,
    with the same degree and lead coefficient, and the span of the roots of $R(x)$
    contain the roots of $Q(x)$, then for any interval $I$ where $R(x)$
    is an obstruction polynomial, $Q(x)$ is also
    an obstruction polynomial, and hence $R(x)$ is not needed.

After construction of these polynomials, we can, for fixed $d, a_d$,
    and $t < a_d^{-1/d}$ find an upper bound for $\Lm(t)$ by
    finding the polynomial $Q$ whose roots
    have the smallest span,  and then appealing to
    Proposition \ref{prop:L-bounds} (a).
This was done in Table \ref{tab:\Lm(t)}, formalized in Theorem
    \ref{thm:\Lm(t)} (a), and displayed in Figure \ref{fig:\Lm(t)}.

Similarly, given this list of polynomials, we can compute the
    least $\ell$ such that any
    interval of length $\ell$ will contain an integer translate of at least
    one of the polynomials in our list.
Then with Proposition \ref{prop:L-bounds} (b) we get an upper bound for
    $\Lp(t)$.
for given $\ell$, we must choose $I_0$ carefully.
If  $I_0$ is too short, we might miss an important obstruction
    polynomial. On the other hand,
if  $I_0$ is long, we will find, along with the obstruction
    polynomials we seek, also (possibly multiple) integer translates of these
    polynomials.
This is inefficient, as we end up doing more
    calculations than we need to.
So we wish to pick  $I_0$ so that it is
    long enough to ensure that we have all important obstruction polynomials,
    and yet small enough that we are not doing more work than necessary.
We do this by ensuring that $I_0$, the interval which contains the roots
    of the polynomials we have found, has the property that $|I_0|$ is just
    greater $\ell+1$.
This ensures that there are no other useful obstruction polynomials
    that we might have missed, since any obstruction polynomial having a span of length $\ell$ will then have some
    integer translate lying in $I_0$.
(We might have to re-run the calculation if $|I_0|$ is too small based on the
    current value of $\ell$.)
We can achieve tighter upper bounds for $\Lp(t)$ by considering the list
of all obstruction polynomials we found such that $a_d^{-1/d} \geq t$.

This computation was done for $t = \frac{1}{2}$ (Table
    \ref{tab:1.48} and Theorem \ref{thm:\Lm(1/2)}) and also
 for $20$ other values of $t$  (Table
\ref{tab:\Lp(t)},  Theorem \ref{thm:\Lm(t)} (b)
    and  Figure \ref{fig:\Lp(t)}).
To save space, the list of relevant polynomials for each $t$
    is not given in the table.
(This information is available upon request from the authors.)

\section{critical polynomials: results and proofs}
\label{sec:crit}

    We first establish a relationship
    between critical polynomials and maximal obstructions. We define
    a {\em maximal} nonmonic critical polynomial of an interval $I$ to be
    a critical polynomial $Q(x)=a_dx^d+\cdots+a_0$
    such that the value $a_d^{-1/d}$ is maximal for $Q$ within the set of
    nonmonic critical polynomials for $I$.
Such a polynomial is well defined, as a result of the following Theorem.

\begin{thm}
\label{thm:critical=maximal}
Suppose that an interval  $I$ has a nonmonic critical polynomial. Then  $I$ has a maximal nonmonic critical polynomial,
 $Q(x) = a_d x^d + \cdots + a_0$ say, and
    furthermore $Q$ is also a maximal obstruction polynomial, so that
     $a_d^{-1/d}$ is the maximal obstruction.
\end{thm}

To prove this result, we will apply the following version of a classical lemma.

\begin{lem}[{\cite[p. 77]{Borwein02}}]\label{lem:was BE}
Let $Q(x)$ and $R(x)$ be two (not necessarily monic) integer polynomials.
Further suppose that $Q(x) = a_d x^d + \cdots + a_0$ is a critical polynomial for the
interval $I$, and that the  integer polynomial $R(x)$
satisfies  $\lh R\rh_I< a_d^{-1/d}$. Then $Q$ divides $R$.
\end{lem}
\begin{proof} From equations (\ref{eqn:resultant}) and
    (\ref{eqn:resultant2}), with $R(x)$ replacing $P(x)$, we must
    have $\Res(Q,R)=0$.
\end{proof}

This result, essentially known to early workers on integer transfinite diameter
(Gor\v skov, Sanov, Trigub, Aparicio Bernardo, ...), has appeared in the literature in various forms -- see for instance
 Chudnovsky \cite[Lemma 2.3]{Chudnovsky83}, Montgomery \cite[Chapter 10]{Montgomery94}, Borwein and
Erd\'elyi \cite{BorweinErdelyi96}, Flammang, Rhin and Smyth \cite{FlammangRhinSmyth97}.

\begin{proof}[Proof of Theorem \ref{thm:critical=maximal}]
We first observe  that nonmonic critical polynomials are
    obstruction polynomials.
Conversely, if an obstruction is greater than $\tz(I)$ then its
associated polynomial is     also a critical polynomial.

Assume that $I$ has a nonmonic critical polynomial, and consider the
nonempty set $\A = \{a_d^{-1/d}\}$ of obstructions coming from the
nonmonic critical polynomials of $I$. Any integer polynomial $R(x)$ (not
necessarily monic), must, by Lemma \ref{lem:was BE}, contain as factors all
critical polynomials $Q$ whose obstructions $a_d^{-1/d}$ are
strictly greater than $\lh R(x)\rh_I$.
Therefore $\lh R(x)\rh_I\geq \ell$ for any limit point $\ell$ of
$\A$, and hence $\tz(I) = \ell$.
So if $\A$ has a limit point, then it must be $\inf(\A)$. Thus $\sup(\A)$ is
attained, and there is a maximal nonmonic
    critical polynomial $Q$ say. Then $Q$ is also a maximal obstruction polynomial.
\end{proof}

\begin{cor}\label{cor:critical=maximal}
Conjecture \ref{conj:maximal} and Conjecture \ref{conj:maximal=tm}
together imply
    Conjecture \ref{conj:critical}.
\end{cor}

\begin{proof} From
the proof above, we see that an obstruction that is greater
    than $\tz(I)$ is associated to a critical polynomial.
The existence of such an obstruction is a consequence of
    Conjecture \ref{conj:maximal} and Conjecture \ref{conj:maximal=tm}.
\end{proof}
\begin{proof}[Proof of Proposition \ref{prop:crit}]  Now $\tz(I)\leq \inf_i a_{d_i,i}^{-{1}/{d_i}}$, by the
 definition of a critical polynomial. But if
this inequality were strict, then we could find an integer
 polynomial $R$ with $\tz(I)\leq \lh R\rh_I<\inf_i a_{d_i,i}^{-{1}/{d_i}}$. But then,
from Lemma  \ref{lem:was BE}, $R$ would have to be divisible by all
the $Q_i$, which is impossible.
\end{proof}

\section{Farey intervals and the proof of Theorem \ref{thm:maximal}}
\label{sec:Farey}

Every closed interval $I$ has a least positive integer $q$ such that some rational $p/q$ with $(p,q)=1$ lies in the interior of $I$.
If $q\ge 2$ then $I$ belongs to a unique Farey interval $\Farey$ whose endpoints are consecutive fractions in the Farey sequence of order $q-1$. We define this interval to be the {\it minimal Farey interval containing $I$}.

 Theorem \ref{thm:maximal} follows directly from our next result.

\begin{thm}
\label{thm:rational} Let $I$ be an interval not containing an
integer in its interior, and $\Farey $
    be the minimal Farey interval containing $I$.
Then $(c_1+c_2) x - (b_1+b_2)$ is a critical polynomial for $I$.
Moreover, the maximal obstruction for $I$ is 
\[ =
\begin{cases}

 \frac{1}{c_1} & \text{ if } c_1\ge 2, \frac{b_1}{c_1}\in I,\frac{b_2}{c_2}\not\in I;\\
 \frac{1}{c_2} & \text{ if } c_2\ge 2, \frac{b_1}{c_1}\not\in I,\frac{b_2}{c_2}\in I;\\
 \frac{1}{\min(c_1,c_2)} & \text{ if } c_1\ge 2, c_2\ge 2 \text{ and } I=\Farey;\\
\frac{1}{c_1+c_2} & \text{ otherwise. }
\end{cases}
\]
\end{thm}

\begin{proof}
Now the polynomial $Q(x) = (c_1 x - b_1)^{c_2} (c_2 x - b_2)^{c_1}$
has a local maximum of
    $\left(\frac{1}{c_1+c_2}\right)^{c_1+c_2}$ at $x= \frac{b_1+ b_2}{c_1+c_2}$ .
     Thus, by
continuity, there exist integers $r_1$ and $r_2$ such that
    $R(x):=Q(x)^{r_1} ((c_1+c_2) x - (b_1+b_2))^{r_2}$ has normalized supremum
    less than $\frac{1}{c_1+c_2}$.
Hence $(c_1+c_2) x - (b_1+b_2)$ is an obstruction polynomial. Now
$\frac{b_1+ b_2}{c_1+c_2}\in I$, as otherwise $I$ would be contained
in one of the Farey intervals $\left[\frac{b_1}{c_1},\frac{b_1+
b_2}{c_1+c_2}\right]$ or $\left[\frac{b_1+
b_2}{c_1+c_2},\frac{b_2}{c_2}\right]$.

Since the polynomials $(c_1+c_2) x - (b_1+b_2)$,  $c_1x-b_1$ and
$c_2x-b_2$ are critical only if their roots
    are in $I$, and are, as factors of $R$, by Lemma \ref{lem:was BE} the
     only three possible maximal critical polynomials in this
    Farey interval, we get the final result.
\end{proof}

\begin{thm}
\label{thm:farey} Let $\Farey$ with $c_1\ge 2$ be a Farey interval, and suppose that
$b_1^2 \equiv \pm 1 \pmod{c_1}$ and $ b_2^2\equiv B \pmod {c_2}$
where
    $c_1^2 |B| < c_2^2$.
Then
$$
\tm\left(\Farey\right) = \frac{1}{c_1}.
$$
\end{thm}

\begin{proof} From
    \cite[p. 1905]{BorweinPinnerPritsker03} we have that there exists a monic quadratic
    integer polynomial $P(x)$ which has the property that
    $P\left(\frac{b_1}{c_1}\right) = \pm\frac{ 1}{c_1^2}$ and
    $P\left(\frac{b_2}{c_2}\right) = \frac{B}{c_2^2}$.
Since its critical point  is at a half integer, it is strictly
monotonic on the Farey interval. Hence it attains its maximum at one
of its endpoints, and $\left|P\left(\frac{b_1}{c_1}\right)\right|>
    \left|P\left(\frac{b_2}{c_2}\right)\right|$.
\end{proof}

\begin{thm}
\label{thm:farey 2} Let $P(x) = x^2 + a_1 x + a_0$ be an irreducible
integer polynomial with
    real roots.
Then there exist infinitely many Farey intervals
    for which $P(x)$ attains the maximal obstruction.
\end{thm}

\begin{proof}
We know (Pell's Equation) that the equation $x^2 + a_1 x y + a_0
y^2 = \pm 1$ has an infinite
    number of solutions
    $(x,y)=(b_i,c_i)$.
These solutions have the property that
    $P\left(\frac{b_i}{c_i}\right)=\pm\frac{ 1}{c_i^2}$. Further, by
    choosing a suitable subsequence we may assume that both the
     $c_i$ and the   $b_i/c_i$ are
    monotonically increasing.
Thus for any interval $I :=
\left[\frac{b_i}{c_i},\frac{b_{i+1}}{c_{i+1}}\right]$
    not containing a half-integer, we see that $P(x)$ attains the maximal
    obstruction $1/c_i$ with $Q(x)=c_ix-b_i$, so that $\tm(I) =1/c_i$.
This happens infinitely often as the $\frac{b_i}{c_i}$ tend to a
root of $P(x)$.

We can find a
$\frac{b}{c}\in\left[\frac{b_i}{c_i},\frac{b_{i+1}}{c_{i+1}}\right]$
   such that $\left[\frac{b_i}{c_i},\frac{b}{c}\right]$ is a Farey interval,
   and hence $P(x)$ attains its
    maximal obstruction $1/c_i$ on this interval.
\end{proof}

It should be noted that this method of proof will not work for
polynomials of degree
     3 or higher, as the resulting Thue equation
    \[ x^n + a_{n-1} x^{n-1} y + \cdots + a_0 y^n = \pm 1\]
    has only a finite number of integer solutions \cite{Sprindzuk93}.

\section{Study of $\tm(b)$}
\label{sec:interval 1/3} In this section we consider intervals
$[0,b]$, with $\tm(b)$ denoting $\tm([0,b])$.
 Our first result for such intervals is a consequence of
 Theorem \ref{thm:rational}.

\begin{cor}
Let $n\ge 2$ and $\frac{1}{n} < b < \frac{1}{n-1}$. Then
$\frac{1}{n}$ is the maximal obstruction of $[0, b]$.
\end{cor}

\begin{thm}
\label{thm:interval 1/n} For all $n \in \BbN$ there exists $\delta_n
>\fracn-\frac{1}{n}$ such that
    for all $0 \leq \varepsilon \leq \delta_n$
    \[\tm\left(\left[0, \tfrac{1}{n} + \varepsilon\right]\right) = \tfrac{1}{n}.\]
\end{thm}

\begin{proof}
Consider the polynomial $P_n(x) = x^{n^2-2} (x^2-n x + 1)$.
It has the following properties:
\begin{itemize}
\item $P_n\left(\frac{1}{n}\right) = \left(\frac{1}{n}\right)^{n^2}$;
\item $P_n(x)$ has a local maximum (with respect to $x$) at $x=\frac{1}{n}$;
\item $P_n(x)$ is strictly increasing (with respect to $x$) on $\left[0,
    \frac{1}{n}\right]$;
\item $P_n(x)$ has a root $\beta_n = \fracn$  strictly
    greater than $\frac{1}{n}$;
\item $P_n(x)$ is strictly decreasing on $\left[\frac{1}{n}, \beta_n\right]$.
\end{itemize}

Let $\alpha_n$ be the minimal root, strictly greater than $\beta_n$, of the
    equation $|P_n(x)| = \frac{1}{n^{n^2}}$.
Thus $P_n(x)$ demonstrates that $\tm( \alpha_n) = \frac{1}{n}$, where
$\alpha_n > \beta_n = 2/(n + \sqrt{n^2-4}) > \frac{1}{n}$.
\end{proof}

\begin{thm}
\label{thm:interval 1/3} We have that
\begin{enumerate}
\item[(a)] $\tm(b) = \tfrac{1}{4}$ for $b \in [\tfrac{1}{4}, 0.303]$;  \label{it:0.3}
\item[(b)] $\tm(b) = \tfrac{1}{3}$ for $b \in [\tfrac{1}{3}, 0.465]$; \label{it:0.46}
\item[(c)] $\tm(b) = \tfrac{1}{2}$ for $b \in [\tfrac{1}{2}, 1.26]$; \label{it:1.26}
\item[(d)] $\tm(1.328) > \tfrac{1}{2}$.  \label{it:1.32}
\end{enumerate}
Hence, in the notation of Theorem \ref{thm:interval 1/n},
$0.76\le\delta_2<0.828$,
$\delta_3 > 0.132  $ and $ \delta_4 > 0.053$.
\end{thm}

\begin{proof}
The optimal monic polynomials needed for Parts
    (a) and (b) are given in
    Table \ref{tab:interval 1/3}. In each case they attain the
    maximal obstruction $\tfrac{1}{4}$ and $\tfrac{1}{3}$ respectively.
As before, a slightly larger interval can be computed exactly,
    by solving $P(x) = \pm \left(\frac{1}{4}\right)^{\deg P}$
    or $P(x) = \pm \left(\frac{1}{3}\right)^{\deg P}$ respectively.
The values of 0.303 and 0.465 have been rounded down to ensure that the
    inequality still holds.
Part (c) follows from the first part of Table \ref{tab:interval 1},
    using the map $x \mapsto 1-x$, with the same comments to the
    exact values as above.
Part (d) is proved using Lemma BPP using the obstruction
    polynomial $7x^3-14x^2+7x-1$.
Here 1.328 is an approximation to its largest root,  rounded up to ensure  that (d) holds.
\end{proof}
The factors used for the construction of the polynomials in Table
    \ref{tab:interval 1/3} were found using the techniques discussed
    in Section \ref{sec:comp}, making use of the constraints
    given by  Lemma \ref{lem:resultant-value}.

\begin{table}
\begin{tabular}{|l|}
\hline
$\tm(b) = \frac{1}{4}$ for $b \in [\tfrac{1}{4}, 0.303]$ by $P(x)=$ \\
$ x^{640}(x^5+432 x^4-456 x^3+179 x^2-31 x+2)^{47}$ \\
$(x^7+8760 x^6-13342 x^5+8388 x^4-2784 x^3+514 x^2-50 x+2)^{35} $ \\
\\

$\tm(b) = \frac{1}{3}$ for $b \in [\tfrac{1}{3}, 0.465]$ by  $P(x)=$\\
$x^{1652706720}(x^7-1233 x^6+2406 x^5-1913 x^4+791 x^3-179 x^2+21 x-1)^{118037088}$ \\
$(x^8+4842 x^7-10935 x^6+10355 x^5-5317 x^4+1594 x^3-278 x^2+26 x-1)
    ^{156479575}$ \\
$(x^8+14184 x^7-34944 x^6+36442 x^5-20832 x^4+7041 x^3-1405 x^2+153 x-7)
   ^{72166388}$ \\
$(x^8+7812 x^7-18072 x^6+17561 x^5-9271 x^4+2864 x^3-516 x^2+50 x-2)
    ^{4378185}$ \\

\hline
\end{tabular}
\caption{Optimal monic integer polynomials used for the proof of Theorem \ref{thm:interval
1/3}.} \label{tab:interval 1/3}
\end{table}


Bounds have been given on the exponents of certain factors for large
    integer Chebyshev polynomials used for estimating $\tz(I)$.
For example, for the interval $I = [0,1]$, Pritsker \cite{Pritsker99}
shows that $(x (1-x))^\gamma$, where $0.2961 \le \gamma \le 0.3634$,
must appear as a factor in any  polynomial $R$ (normalized to have
degree $1$), for which $\lh R\rh_I$ is sufficiently close to $\tz(I)$.

Following \cite{FlammangRhinSmyth97}, we now determine a lower bound
for $\gamma(b)$ such that $x^{\gamma(b)}$ must divide any normalized
monic integer  polynomial $P$ such that $\lh P\rh_{[0,b]}$
approximates $\tm(b)$ sufficiently closely.

Suppose that the function $m(b)$ is an upper bound for $\tm(b)$. Then
by Proposition 5.3 and Lemma 5.2 of \cite{FlammangRhinSmyth97}
    we have that $\gamma(b)$ is bounded below by the least positive root of
\[ \frac{(1+x)^{1+x}}{(1-x)^{1-x} (2x)^{2x} b^x} = \frac{1}{m(b)}. \]

\begin{figure}[ht]
\[\psfig{file= 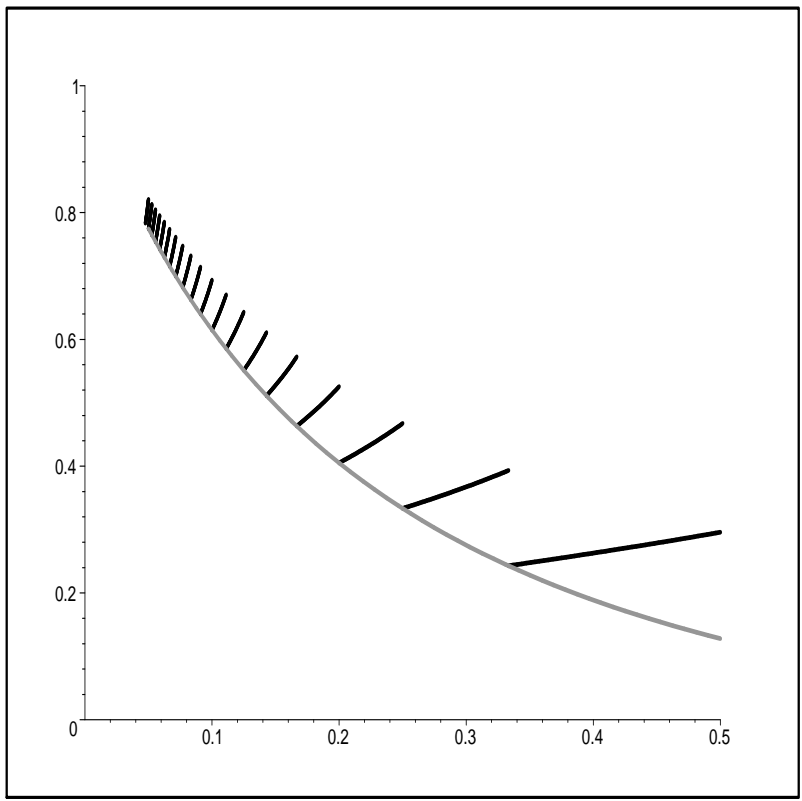,height= 300pt,width= 300pt,angle=270}\]
\caption{Lower bounds for $\gamma(b)$. \newline grey line -- lower
bound using $m(b) = b$ (known), \newline solid line -- lower bound
assuming $m(b) = \frac{1}{\lceil 1/b \rceil}$ (Conjecture \ref{conj:[0,b]}).} \label{fig:gamma}
\end{figure}

So in particular, if $\tm(b) = \frac{1}{\lceil 1/b \rceil}$ for
$b\in[0,1]$ as in
    Conjecture \ref{conj:[0,b]}, then our lower bound for $\gamma(b)$ would
    have infinitely many discontinuities in this range
    (Figure \ref{fig:gamma} -- black lines).
However, we know, by using the polynomial $x$, that we have a provable,
    albeit weaker, upper bound $m(b)=b$.
This gives us a proven lower bound for $\gamma(b)$
    (Figure \ref{fig:gamma} -- grey line).

\begin{prop}
\label{thm:gamma(b)} We have $\displaystyle \lim_{b\to 0} \gamma(b)
= 1$.
\end{prop}

\begin{proof}
Define \[ T(x,b) = \frac{(1+x)^{1+x}}{(1-x)^{1-x} (2x)^{2x} b^x} -
                      \frac{1}{b}. \]

Now $T(x,b)$ has a positive local maximum  at
   $x = \frac{1}{\sqrt{1+4 b}} \to 1$ as $b \to 0$,while
   $T(1-\sqrt{b}, b) < 0$ for $0< b < 0.04$, so that $T(x,b)=0$ has a root 
   in $[1-\sqrt{b},\frac{1}{\sqrt{1+4 b}}]$. Further, since $T(x,b)$ is increasing for $x\in[0,\frac{1}{\sqrt{1+4 b}}]$ this root is the least positive root of $T(x,b)=0$. Hence     
    $\gamma(b)>1-\sqrt{b}$,
giving the result.
\end{proof}


\section{Proof of Counterexample \ref{thm:unattainable}}
\label{sec:unattainable}

For the proof of Counterexample  \ref{thm:unattainable}
    we need the following $p$-adic result.

\begin{prop} \label{prop:p-adic}Suppose that $Q(x)=a_d x^d + \cdots +a_0 \in\BbZ[x]$ is
a maximal obstruction polynomial for the interval $I$, and that the
maximal obstruction is attained by some monic integer polynomial
$P(x)$. Then $\gcd(a_0,a_d)=1$ and, for every prime $p$ dividing $a_d$
we have
$$
\left|\frac{a_{d-i}}{a_d}\right|_p\le\left|\frac{1}{a_d}\right|_p^{i/d}
\quad (i=0,\cdots,d).
$$

In particular, if $a_d$ is square-free then $\frac{1}{a_d}(Q(x)-Q(0))$ has
integer coefficients.

\end{prop}
Here $|.|_p$ is the usual $p$-adic valuation on $\BbQ$. For the proof, it is extended to $\overline\BbQ$.

\begin{proof}
Take $\beta$ to be any root of $Q(x)$, and $p$ any prime factor of
$a_d$. Let $P(x)$ be of degree $m$. Then, as $P(x)$ attains the
obstruction, $P(\beta)=\pm a_d^{-m/d}$, so that
$|P(\beta)|_p=|1/a_d|_p^{m/d}>1$.   If $|\beta|_p\le 1$ then
$|P(\beta)|_p\le 1$, a contradiction, as $P(x)$ has integer
coefficients. Hence $|\beta|_p>1$ and
$|P(\beta)|_p=|\beta|_p^m=|1/a_d|_p^{m/d}$, giving
\begin{equation}\label{E-1}
|\beta|_p=|a_d|_p^{-1/d}.
\end{equation}
Applying (\ref{E-1}) for all roots $\beta_j$ of $Q(x)$ we get
$|\prod_j\beta_j|_p=|1/a_d|_p$. But also from
$a_d^{-1}Q(x)=\prod_j(x-\beta_j)$ we have that
$|\prod_j\beta_j|_p=|a_0/a_d|_p$. Hence $|a_0|_p=1$. Doing this for
all $p|a_d$ we obtain $(a_0,a_d)=1$.
Furthermore, if
$\left|\frac{a_{d-i}}{a_d}\right|_p >
\left|\frac{1}{a_d}\right|_p^{i/d}$ for any $i$ then the Newton
polygon of $P$ (see for instance \cite[p. 73]{Weiss63}) tells us that $|\beta_j|_p > |1/a_d|_p^{1/d}$ for
some $j$, contradicting (\ref{E-1}).

In the case of $a_d$ square-free,
$\left|\frac{1}{a_d}\right|_p^{i/d}< p$ for $1\le i<d$, so that
$\left|\frac{a_{d-i}}{a_d}\right|_p\le 1$, and hence, using all
primes $p$ dividing $a_d$, we see that $\frac{a_{d-i}}{a_d}$ is an integer.
\end{proof}

\begin{proof}[Proof of Counterexample \ref{thm:unattainable}]

The fact that $7 x^3 + 4 x^2 - 2 x -1$ is a maximal obstruction
polynomial for
    the interval $I=[-0.684, 0.517]$ can be verified by showing that
    it is a critical polynomial.
This follows from the fact that the  polynomial
\begin{eqnarray*}
R(x)&=&  x^{28728} (5 x^3+4 x^2-x-1)^{3739} (7 x^3+4 x^2-2 x-1)^{1140} \\
&& (x^6-24 x^5-20 x^4+10 x^3+9 x^2-x-1)^{420}\\
&& (3 x^5+16 x^4+3 x^3-8 x^2-x+1)^{399}
\end{eqnarray*}
has $\lh R\rh_I<7^{-1/3}$, so that $\tz(I) < 7^{-1/3}$. As
    $7 x^3 + 4 x^2 - 2 x - 1$ has all its roots in $I$, it is therefore
    a critical polynomial.
As always, the interval is an approximation only, and a tighter
   one can easily be computed.

    We now claim that $7 x^3+4 x^2-2 x-1$ is the maximal nonmonic critical polynomial for
    $I$. For any critical polynomial $a_dx^d+\cdots+a_0$ for $I$ with
    $a_d^{-1/d}>\lh R\rh_I$ must be a factor of $R$, by Lemma \ref{lem:was BE}.
But among the four irreducible
    factors of $R$, $7 x^3+4 x^2-2 x-1$ is the only one having all
    its roots within $I$. As it is nonmonic, it must indeed be
the maximal nonmonic critical polynomial for
    $I$. By Theorem \ref{thm:critical=maximal}, this polynomial is
    the maximal obstruction polynomial. However, $\frac{1}{7}(7 x^3+4 x^2-2
    x)$ does not have integer coefficients so that, by
    Proposition \ref{prop:p-adic}, the interval has no optimal monic integer Chebyshev polynomial.
    \end{proof}

\section{Some Final Comments on the Computations and Figures}

Consider Figure \ref{fig:\Lm(t)}.
We see that $\Lm(t) = 0$ for for $t < \frac{1}{2}$, and further that
    $\Lm(t) = 4 t$ for $t > 1$.
So in fact the area of interest is for $t$
    between $\frac{1}{2}$ and $1$.
That being said, the upper bound is only given up to approximately
    $0.89$.
This is because the upper bound from Proposition \ref{prop:L-bounds}(a) is given by high degree polynomials with small
    lead coefficient.
In our search, we compute only up to degree 6.
As $2^{-1/6} \approx 0.89$ this is the limit to our knowledge of the upper
    bound.
If we wished to extend these calculations, we could extend the knowledge
    of the upper bound, but the computation time becomes excessive.
For example, even if we computed up to degree 10, which is
    probably beyond our computational range, we would only get up to
    $0.933$.
As it was, the computations up to degree 6  took over 3000
    CPU hours, and the computation time approximately triples for
    each additional degree.
Similar comments apply to bounding $\Lp(t)$ (Figure \ref{fig:\Lp(t)}) for $t$ close
    to $1$.
In this case, it actually turned out that none of the polynomials with
    lead coefficient 2 and degree 6 were useful in the calculations
    for such $t$, and hence we only get an upper bound for $L_+(t)$ for $t$ up to   $t=2^{-1/5} \approx 0.871$.

While we know from Lemma \ref{lem:LL}(c) that $\Lp(t)\geq 2t$ for $t
\le 1/2$, we do not know $\Lp(t)$ exactly in this range. In order to
get an upper bound for $\Lp(t)$ in at least part of this range, it
would in principle be possible to extend the calculation downwards
from $t=\tfrac{1}{2}$. The lower bound of $\frac{1}{2}$ for $t$ was
chosen, as we computed obstruction
    polynomials of degree $d$, with coefficients up to $2^d$.
If we were to compute up to $3^d$ instead, we would be able to extend this
    graph down to $t=\frac{1}{3}$.
This would, however, be a massive undertaking,
because  we would
    have  $3^6/2^6>11$  times as many possible lead coefficients.
Furthermore, we observed that, for a given degree, the computations took longer
     the higher the lead coefficient was, so this
    factor  $11$ is probably an underestimate.

It may be possible to extend these calculations though in a more sophisticated
    manner, somehow doing a less extensive and more intelligent search for
    obstruction polynomials of higher degree or larger lead
    coefficients.
This would be a worthwhile project, and could lead to some interesting
    new results.

Lastly, consider Figure \ref{fig:gamma}.
This could very easily have been extended all the way to 0.
The reason that we chose not to do this is because the
    hypothetical lower bound (the black lines) starts to merge into itself,
    and the Figure becomes unreadable.
(The lower bound jumps at every $\frac{1}{n}$ which get more frequent as
    $\frac{1}{n} \to 0$.)

{\em Acknowledgement.} We thank the referee for helpful comments.

\clearpage

\providecommand{\bysame}{\leavevmode\hbox to3em{\hrulefill}\thinspace}
\providecommand{\MR}{\relax\ifhmode\unskip\space\fi MR }

\providecommand{\MRhref}[2]{%
  \href{http://www.ams.org/mathscinet-getitem?mr=#1}{#2}
}
\providecommand{\href}[2]{#2}

\end{document}